\documentclass[11pt, letterpaper, leqno]{amsart} 

\usepackage{graphicx} 
\usepackage{amsmath,amsthm,amssymb} 
\usepackage{amsfonts} 
\usepackage{dsfont}
\usepackage{mathrsfs} 
\usepackage[T1]{fontenc} 
\usepackage[utf8]{inputenc} 

\usepackage[dvipsnames, svgnames]{xcolor} 
\usepackage[shortlabels]{enumitem} 
\usepackage[colorlinks,citecolor=citec,hypertexnames=false,urlcolor=urlc,breaklinks=true, pagebackref]{hyperref} 

\usepackage{accents} 

\usepackage{amsbsy} 
\usepackage{amscd} 
\usepackage{latexsym} 
\usepackage{txfonts} 

\usepackage{exscale} 
\usepackage{bbm} 
\usepackage{mathtools} 
\usepackage{bbding}

\usepackage{geometry} 

\makeatletter
\providecommand\@dotsep{4.5}
\makeatother

\usepackage{thmtools} 

\newif\ifsoldark
\newif\ifsollight
\newif\ifclassic
\newif\ifplain

\usepackage{cite} 

\renewcommand*\backref[1]{\ifx#1\relax \else (p. #1) \fi} 

\numberwithin{equation}{section} 

\newcounter{questions}

\theoremstyle{plain}
\newtheorem{theorem}[equation]{Theorem}
\newtheorem{question}[questions]{Question}
\newtheorem{lemma}[equation]{Lemma}
\newtheorem{corollary}[equation]{Corollary}
\newtheorem{proposition}[equation]{Proposition}

\theoremstyle{definition}
\newtheorem{definition}[equation]{Definition}

\theoremstyle{remark}
\newtheorem{remark}[equation]{Remark}

\newcommand{\eps}{\varepsilon}

\newcommand{\dint}{\int\!\!\!\!\!\int}
\newcommand{\dist}{\operatorname{dist}}

\newcommand{\dv}{\operatorname{div}}
\newcommand{\Di}{\operatorname{(D)}}
\newcommand{\n}[1]{\mathscr{#1}}
\newcommand{\bb}[1]{\mathbb{#1}}

\newcommand{\medcup}{\textstyle\bigcup}
\newcommand{\medsum}{\textstyle\sum}

\newcommand{\RNum}[1]{\uppercase\expandafter{\romannumeral #1\relax}}

\DeclareMathOperator{\supp}{supp}

\DeclareMathOperator{\diam}{diam}

\def\div{\mathop{\operatorname{div}}\nolimits}

\def\Xint#1{\mathchoice
	{\XXint\displaystyle\textstyle{#1}}%
	{\XXint\textstyle\scriptstyle{#1}}%
	{\XXint\scriptstyle\scriptscriptstyle{#1}}%
	{\XXint\scriptscriptstyle\scriptscriptstyle{#1}}%
	\!\int}
\def\XXint#1#2#3{{\setbox0=\hbox{$#1{#2#3}{\int}$}
		\vcenter{\hbox{$#2#3$}}\kern-.5\wd0}}

\def\dashint{\Xint-}

\def\Yint#1{\mathchoice
	{\YYint\displaystyle\textstyle{#1}}%
	{\YYint\textstyle\scriptstyle{#1}}%
	{\YYint\scriptstyle\scriptscriptstyle{#1}}%
	{\YYint\scriptscriptstyle\scriptscriptstyle{#1}}%
	\!\dint}
\def\YYint#1#2#3{{\setbox0=\hbox{$#1{#2#3}{\iint}$}
		\vcenter{\hbox{$#2#3$}}\kern-.51\wd0}}
\def\longdash{\mkern-3.5mu{-}\mkern-3.5mu{-}} 

\def\fiint{\Yint\longdash}

\newcommand{\ra}{\rightarrow}
\newcommand{\lra}{\longrightarrow}
\newcommand{\m}[1]{\mathcal{#1}}

\newcommand{\vertiii}[1]{{\left\vert\kern-0.25ex\left\vert\kern-0.25ex\left\vert #1 
		\right\vert\kern-0.25ex\right\vert\kern-0.25ex\right\vert}}
\newcommand{\ADR}{\operatorname{ADR}}

\newcommand{\loc}{\operatorname{loc}}

\newcommand{\dyadic}{\textnormal{dyadic}}

\allowdisplaybreaks

\setlist{nosep} 

\colorlet{citec}{blue}
\colorlet{urlc}{blue}
\colorlet{toc}{blue}
\colorlet{hyperc}{blue}
\colorlet{bpcolor}{NavyBlue}
\colorlet{smcolor}{purple}
\colorlet{impcolor}{blue}
\colorlet{eqcolor}{black}
\colorlet{lmcolor}{black}
\colorlet{propcolor}{black}
\colorlet{thmcolor}{black}
\colorlet{defcolor}{black}
\colorlet{rmcolor}{black}
\colorlet{excolor}{black}

\usepackage{everypage} 
\usepackage{tikzpagenodes} 

\newcommand{\ms}{\medskip}

\newcommand{\R}{\mathbb{R}}

\newcommand{\B}{\mathcal{B}}
\newcommand{\T}{\mathcal{T}}

\newcommand{\sm}{\setminus}

\newcommand{\wt}{\widetilde}
\newcommand{\wh}{\widehat}

\newcommand{\1}{{\mathds 1}}

\newcommand{\dr}{\partial}

\newcommand{\A}{\mathcal{A}}
\newcommand{\D}{\mathbb{D}}

\newcommand{\bp}{\noindent {\it Proof}.\,\,}
\newcommand{\ep}{\hfill$\Box$ \vskip 0.08in}

\begin{document}

\author[J. Feneuil]{Joseph Feneuil}
\address{Laboratoire de math\'ematiques d'Orsay, 
	\\
	Universit\'e Paris-Saclay, CNRS, 
	\\
	91405, Orsay, France} 
\email{joseph.feneuil@gmail.com}

\author[B. Poggi]{Bruno Poggi}
\address{Departament de  Matemàtiques, Universitat Autònoma de Barcelona, Bellaterra,  Catalonia} 
\email{bgpoggi.math@gmail.com}

\thanks{J.F.\   was partially supported by the Simons Foundation grant 601941, GD. B.P.\ was supported by the University of Minnesota Doctoral Dissertation Fellowship grant and in part by the Simons Collaborations in MPS 563916, SM. The authors would like to thank Svitlana Mayboroda and Max Engelstein for insightful conversations.}

\title[Generalized Carleson perturbations of elliptic operators]{Generalized Carleson perturbations of elliptic operators and applications}
\date{\today}

\begin{abstract}  We extend in two directions the notion of perturbations of Carleson type for the Dirichlet problem associated to an elliptic   real second-order divergence-form (possibly degenerate, not necessarily symmetric) elliptic operator.
First, in addition to the classical perturbations of Carleson type, that we call additive Carleson perturbations, we introduce scalar-multiplicative and antisymmetric Carleson perturbations, which both allow non-trivial differences at the boundary.
Second, we consider domains which admit an elliptic PDE in a broad sense: we count as examples the 1-sided NTA (a.k.a. uniform) domains satisfying the capacity density condition, the 1-sided chord-arc domains, the domains with low-dimensional Ahlfors-David regular boundaries, and certain domains with mixed-dimensional boundaries; thus our methods provide a unified perspective on the Carleson perturbation theory of elliptic operators.

Our proofs do not   introduce sawtooth domains or the extrapolation method. We also present several applications to some Dahlberg-Kenig-Pipher operators, free-boundary problems, and we provide a new characterization of $A_{\infty}$ among  elliptic measures. 
\end{abstract}

\maketitle

\ms\noindent{\bf Key words.}
elliptic measures, $A_\infty$-absolute continuity, Carleson perturbations, degenerate operators.

\ms\noindent


{
	\hypersetup{linkcolor=blue}
	\tableofcontents
}
\hypersetup{linkcolor=blue}

\section{Introduction}

In this article, we study additive, scalar-multiplicative, and antisymmetric perturbations of Carleson type for the Dirichlet problem for real second-order divergence-form (possibly degenerate, not necessarily symmetric) elliptic operators on domains which admit an elliptic PDE theory. We call such domains \emph{PDE friendly} (see Section \ref{S2} for our axioms and examples of PDE friendly domains). Roughly speaking, if $L_0$ and $L_1$ are two elliptic operators on such a domain, we seek conditions on the relative structure of $L_1$ to $L_0$ that preserve certain ``good estimates'' for the Dirichlet problem.  In particular, we develop Carleson  perturbations which allow for non-trivial differences at the boundary.

Before describing our full results, which are of a somewhat general nature, let us first review the well-understood situation of the half-plane (but even in this case, some of our results are new),   and the relevant history of Carleson-type perturbations. The reader may also skip directly to Section \ref{sec.results} for our results.

\subsection{A brief review of the Dirichlet problem in the half-plane} \subsubsection{The two Dirichlet problems: continuous data or rough data} Throughout this section, we fix $\Omega=\bb R^n_+=\{(x,t):x\in\bb R^{n-1}, t\in(0,\infty)\}$, $n\geq2$, and we let $A$ be an $n\times n$ matrix of real measurable coefficients on $\Omega$ satisfying the following uniform ellipticity and boundedness conditions
\begin{equation}\label{eq.elliptic}
\tfrac1{C_L}|\xi|^2\leq A(X)\xi\cdot\xi,\qquad |A(X)\xi\cdot\zeta|\leq C_L|\xi||\zeta|,\qquad\text{for each }\xi,\zeta\in\bb R^n,~X\in\Omega.
\end{equation}
Given a matrix $A$, a second-order divergence-form elliptic operator $L$ on $\Omega$ is formally defined as $L=-\dv A\nabla$, and the equation $Lu=0$ in $\Omega$ is interpreted in the weak sense for $u\in W^{1,2}_{\loc}(\Omega)$. Associated to each elliptic operator $L$, there exists a family of Borel probability measures $\{\omega^X_L\}_{X\in\Omega}$ on $\partial\Omega=\bb R^{n-1}\times\{0\}$ so that for any compactly supported continuous function $f$ on $\partial\Omega$, the solution $u\in W^{1,2}_{\loc}(\Omega)\cap C(\overline\Omega)$ to the Dirichlet problem
\begin{equation}\label{eq.dirichletcont}
\left\{\begin{matrix}Lu=0,\qquad\text{in }\Omega,\\u=f,\qquad\text{on }\partial\Omega,\end{matrix}\right.
\end{equation}
may be written as
\begin{equation}\label{eq.represent}\nonumber
u(X)=\int_{\partial\Omega}f\,d\omega_L^X,\qquad\text{for each }X\in\Omega.
\end{equation}
The measure $\omega^X$ on $\partial\Omega$ is called the \emph{elliptic measure} of $\Omega$ associated with the operator $L$ and with pole at $X$. For the half-plane, if $A\equiv I_{n\times n}$ so that $L=-\Delta$, then the $L-$elliptic measure, known in this special case as the \emph{harmonic measure}, is mutually absolutely continuous with respect to the surface measure $\sigma=\n L^{n-1}$ on the boundary, and the Radon-Nikodym derivative $k^X_{-\Delta}= d\omega_{-\Delta}^X/d\sigma$, known as the \emph{Poisson kernel}, satisfies certain scale-invariant reverse H\"older inequalities. We write $\omega_{-\Delta}\in A_{\infty}(\sigma)$ to denote this quantitative absolute continuity property (and see Section \ref{S2} for precise definitions).  These reverse H\"older inequalities allow one to solve the Dirichlet problem with rough data on the boundary. More precisely, if $L=-\Delta$, then for each $p\in(1,\infty)$ and each $f\in L^p(\partial\Omega)$, there exists $u\in W^{1,2}_{\loc}(\Omega)$ such that
\begin{equation}\label{eq.dirichletp}
\left\{\begin{matrix}Lu=0,&\qquad\text{in }\Omega,\\u\lra f,&\qquad\text{non-tangentially }\sigma-\text{a.e.},\\ \Vert N(u)\Vert_{L^p(\partial\Omega,\sigma)}\leq C\Vert f\Vert_{L^p(\partial\Omega,\sigma)},\end{matrix}\right.
\end{equation}
where $N(u)(x):=\sup_{Y\in\gamma(x)}|u(Y)|$, $~\gamma(x):=\big\{Y\in\Omega: |Y-x|\leq2\delta(Y)\big\}$, for $x\in\partial\Omega$, and
\begin{equation}\label{eq.delta}
\delta(Y):=\dist(Y,\partial\Omega),\qquad Y\in\Omega,
\end{equation}
and $u\ra f$ non-tangentially at $x\in\partial\Omega$ if $\lim_{\gamma(x)\ni Y\ra x}u(Y)=f(x)$. The function $N(u)$ is known as the \emph{non-tangential maximal function} of $u$. 

Note that the continuous Dirichlet problem (\ref{eq.dirichletcont}) on the half-plane is solvable for any elliptic operator $L$ whose matrix satisfies (\ref{eq.elliptic}). Naturally, one may wonder whether, for some $p>1$, the  Dirichlet problem with $L^p$ data (\ref{eq.dirichletp}) (henceforth referred to as $\Di_p$) is solvable for elliptic operators other than the Laplacian. It turns out that the question of $\Di_p$ solvability for an operator $L$ is equivalent to whether   $\omega_L\ll\sigma$ and $k^X_L=d\omega_L^X/d\sigma$ satisfies a scale invariant $p'-$reverse H\"older inequality, where $\frac1p+\frac1{p'}=1$, and therefore the machinery of the elliptic measure is a sensible means  to attack the Dirichlet problem. In particular, if one could find an elliptic operator $L$ whose elliptic measure is singular with respect to the surface measure, then it would follow that for this $L$, one cannot solve $\Di_p$ for any $p>1$. The existence of such an $L$  is precisely the pivotal result of Caffarelli-Fabes-Kenig \cite{cfk} (via the Beurling-Ahlfors theory on quasiconformal mappings), and independently, Modica-Mortola \cite{mm}.

\subsubsection{Conditions that guarantee the absolute continuity of elliptic measure with respect to surface measure} The aforementioned examples show that we must place conditions on the matrix $A$ to guarantee the solvability of $\Di_p$ for some $p\in(1,\infty)$ (or, equivalently, that $\omega_L\in A_{\infty}(\sigma)$). The conditions that have historically been considered for real elliptic operators can roughly be categorized into either $t$-independent,    regularity,  or   perturbative assumptions. In this paper, we are interested mainly in the latter, but let us say a few words about the former two. 

The $t-$independence assumption is a natural starting place owing to  \cite{cfk}, where they show that a square Dini condition on the transversal modulus of continuity of $A$ is necessary in order to have solvability of the Dirichlet problem with rough data (a few years later, Fabes-Jerison-Kenig \cite{fjk} obtained the sufficiency of this condition). Moreover, this is the situation that arises from the pullback of the Laplacian on a domain above a Lipschitz graph via the mapping that ``flattens'' the boundary.  The problem $\Di_2$ for $t-$independent real symmetric matrices was solved by Jerison and Kenig in \cite{jk1} (grounded in the pioneering work of Dahlberg \cite{dah1, dah2} for the Laplacian on Lipschitz domains). Later, for the $t-$independent real non-symmetric matrices, $\Di_p$ for sufficiently large $p$ has been solved by Kenig-Koch-Pipher-Toro \cite{kkpt} ($n=2$) and Hofmann-Kenig-Mayboroda-Pipher \cite{hkmps} ($n\geq3$).  We also note that, as pioneered by Fabes-Jerison-Kenig \cite{fjk}, a lot of work for $\Di_2$ has been done in the case of complex-valued operators with $t-$independent coefficients, essentially by perturbing (in $L^{\infty}$ norms) from  the real  case. Usual techniques have been either by the method of layer potentials \cite{aaahk, bhlmp, bhlmp2}, or a functional calculus of Dirac-type operators \cite{aah, aamc}, and these are strong enough to also yield solvability results for the Neumann and Regularity problems; see   \cite{bhlmp} for further discussion.

The regularity condition is borne out from a conjecture posed by Dahlberg in 1984. Dahlberg, Kenig and Stein constructed \cite{dah4} a one-to-one mapping from a Lipschitz domain onto the half-plane for which the pullback of the Laplacian results in a symmetric elliptic operator $L=-\dv A\nabla$ on the half-plane $\Omega$ satisfying (recall $\delta$ is defined in (\ref{eq.delta}))
\begin{enumerate}[({A}1)]
	\item\label{item.linfty} $\delta\nabla A\in L^{\infty}(\Omega)$, and
	\item\label{item.cm} $\delta|\nabla A|^2\,d\n L^n$ is a Carleson measure on $\Omega$; that is, there exists   $C>0$ so that for each $x\in\partial\Omega$ and $r>0$, if $B(x,r)$ is a ball in $\bb R^{n}$  , we have that
	\[
	\dint_{B(x,r)}\delta(Y)|\nabla A(Y)|^2\,dY\leq Cr^{n-1}.
	\]
\end{enumerate}
Since Dahlberg had shown in his celebrated work \cite{dah1, dah2} that $\Di_2$ was solvable for the Laplacian on a Lipschitz domain, he reasonably conjectured that  $\Di_2$ is solvable for any real symmetric elliptic matrix $A$ satisfying the assumptions \ref{item.linfty}-\ref{item.cm}. This question would be resolved over a decade later by Kenig and Pipher \cite{kp3}, and   the real elliptic operators whose matrices satisfy \ref{item.linfty}-\ref{item.cm} have since come to be known as the  Dahlberg-Kenig-Pipher (DKP)  operators. These regularity assumptions are close to optimal (see \cite[Theorem 4.11]{fkp}, \cite{pog1}, and \cite{hmmtz}), but we will revisit these considerations for certain DKP operators further below.  Lastly, we do mention that, by assuming some smallness of the Carleson measure in \ref{item.cm}, Dindos-Petermichl-Pipher \cite{dpp} have obtained the solvability of $\Di_p$ for $p\in(1,\infty)$.

Other than the $t-$independent and regularity conditions, it is natural to  ask whether  the absolute continuity property should be stable under some perturbations of the matrices, although this raises the question of  what type of perturbation to consider. Let us be   more precise: suppose that $L_0$ and $L$ are two elliptic second-order divergence form operators on $\Omega$, with associated matrices $A_0$ and $A$, and elliptic measures $\omega_0$ and $\omega$, respectively. 

\begin{question}\label{q1} What conditions may we ask of the pair $(A,A_0)$  so that if $\omega_0\in A_{\infty}(\sigma)$, then $\omega\in A_{\infty}(\sigma)$?
\end{question}

This question is our main object of study (which we will consider in the generality of PDE friendly domains), so let us now review its history; for a similar and excellent review, see \cite{ahmt}. 

\subsection{History of the Carleson perturbations} \subsubsection{Early results} The first results in this direction are found in \cite{fjk, dah3}. In the setting where $\Omega$ is the unit ball  in $\bb R^n$, the condition that Dahlberg asked of the pair $(A,A_0)$ of symmetric  operators is that the \emph{disagreement} $\rho(A,A_0)$ defined as
\begin{equation}\label{eq.disagreement}\nonumber
\rho(A,A_0)(X):=\sup_{Y\in B(X, \delta(X)/2)}|A(Y)-A_0(Y)|,\qquad  X\in\Omega,
\end{equation} 
satisfies the following \emph{vanishing Carleson measure  condition}
\begin{equation}\label{eq.vanishing}
\lim_{r\searrow0}\sup_{x\in\partial\Omega}h(x,r)=\lim_{r\searrow0}\sup_{x\in\partial\Omega}\Big(\frac1{\sigma(B(x,r)\cap\partial\Omega)}\dint_{B(x,r)\cap\Omega}\frac{\rho(A,A_0)^2(X)}{\delta(X)}\,dX\Big)^{\frac12}=0,
\end{equation}
where $\sigma$ is the Hausdorff $(n-1)-$dimensional measure on the unit sphere $\partial\Omega$. In this case, if $\omega_0\in A_{\infty}(\sigma)$ and its Poisson kernel $k_0=d\omega_0/d\sigma\in RH_p$ (see Proposition \ref{pr.Ainfty}), then $\omega\ll\sigma$ and $k=d\omega/d\sigma\in RH_p$, so that the solvability of $\Di_{p'}$   is stable (with the same $p'$) under the condition (\ref{eq.vanishing}). The fact that the reverse H\"older exponent is preserved by (\ref{eq.vanishing}) suggests that there might be a weaker condition than (\ref{eq.vanishing}) which preserves the $A_{\infty}$ membership but not the $RH$ exponent. Fefferman \cite{fe89} thus showed a few years later that, again in the context of symmetric operators on the unit ball, if $\omega_0\in A_{\infty}(\sigma)$, and if the \emph{area integral}
\begin{equation}\label{eq.defareaint}\nonumber
\n A(\rho(A,A_0))(x):=\Big(\dint_{\gamma(x)}(\rho(A,A_0)(X))^2\frac{dX}{|B(X,\delta(X)/2)|}\Big)^{1/2},\qquad x\in\partial\Omega,
\end{equation} 
satisfies
\begin{equation}\label{eq.areaint}
\n A(\rho(A,A_0))\in L^{\infty}(\partial\Omega,\sigma),
\end{equation}
then $\omega\in A_{\infty}(\sigma)$. It is clear that (\ref{eq.areaint}) is not a vanishing condition; moreover, via Fubini's theorem, one can see that (in the case of the unit ball)
\begin{multline}\label{eq.fubini}\nonumber
h(x,r)\lesssim\Big(\frac1{\sigma(B(x,Cr))\cap\partial\Omega)}\dint_{B(x,Cr)\cap\partial\Omega}\n A(\rho(A,A_0))(x)^2\,d\sigma\Big)^{1/2}\\ \leq\Vert\n A(\rho(A,A_0))\Vert_{L^{\infty}(\partial\Omega,\sigma)},
\end{multline}
and it would be shown in \cite{fkp} that (\ref{eq.areaint}) does not preserve the $RH$ exponent. Next, one may wonder whether (\ref{eq.areaint}) is an optimal condition on $\rho(A,A_0)$ that guarantees the stability of the $A_{\infty}$ property. But the answer to this question is \emph{no}: Fefferman-Kenig-Pipher \cite{fkp} showed that the optimal assumption (at least in the cases of the unit ball or half-plane) which preserves the $A_{\infty}$ property is  that $\rho(A,A_0)^2\delta^{-1}$ is the density of a Carleson measure; in other words, the optimal condition is that
\begin{equation}\label{eq.hcm}
\sup_{r\in(0,\diam(\partial\Omega))}\sup_{x\in\partial\Omega}h(x,r)<+\infty.
\end{equation}
With the landmark paper  of \cite{fkp}, one could say that the contemporary era of the perturbation results was launched: since then, the perturbation results have often assumed  variants of the Carleson measure hypothesis (\ref{eq.hcm}).

\subsubsection{Optimality of the FKP condition on the disagreement} Their proof of optimality  relied on a newfound characterization of $A_{\infty}$ on $\bb R^n$ via a Carleson measure property. Their characterization can be formulated as follows \cite[Page 225]{ste}: Suppose that $w\in L^1_{\loc}(\bb R^n)$ is a non-negative function such that the measure $w\,dx$ is doubling on $\bb R^n$, and that $\Phi$ is a non-negative Schwartz function with $\int_{\bb R^n}\Phi\,dx=1$. Then $w\in A_{\infty}$ if and only if  
\begin{equation}\label{eq.fkpchar}
d\mu:=\frac{|\nabla_x(w*\Phi_t)|^2}{|w*\Phi_t|^2}t\,dx\,dt
\end{equation}
is  a Carleson measure in $\bb R^{n+1}_+$ (here,  $\partial\bb R^{n+1}_+$ is endowed with the $n-$dimensional Hausdorff measure). They go on to show   characterizations of $A_p$ and $RH_p$ through similar Carleson measure conditions \cite[Theorem 3.3]{fkp}. This characterization of $A_{\infty}$ via Carleson measures is not \emph{too} surprising, owing to the classical results that $w\in A_{\infty}$ implies $\log w\in BMO$, and the square-function characterization of $BMO$ \cite{ste}. Their result also fits as a multiplicative analogue of the classical theory of differentiation and a condition of Zygmund; see \cite{fkp} for further discussion. Still, we note that this characterization is for the ``classical'' $A_{\infty}$ space of non-negative weights on Euclidean space, and thus we entertain  

\begin{question}\label{q2} Are there characterizations of $A_{\infty}$ among doubling measures on rough boundaries of domains via Carleson measure properties?
\end{question}

Furthermore, we remark that the FKP condition (\ref{eq.hcm}) is optimal \emph{as a condition on the disagreement function $\rho(A,A_0)$}, which is a scale-invariant version of the difference of the two matrices. In fact, two matrices $A$ and $A_0$ for which $\rho(A,A_0)$ satisfies (\ref{eq.hcm}) must necessarily agree almost everywhere at the boundary of the domain. Therefore, the Carleson perturbations of \cite{fkp} are not adequate to deal with non-trivial perturbations at the boundary. This observation raises 

\begin{question}\label{q3} Could there be a different type of perturbation which allows for a non-trivial difference of the matrices at the boundary?
\end{question}

Note that if $A_0$ is $t-$independent, then the result of \cite{kkpt} guarantees that \emph{any} $t-$independent perturbation from $A_0$ which maintains the ellipticity conditions (\ref{eq.elliptic}) will also preserve the $A_{\infty}$ property; and so our question admits a well-known positive answer among the $t-$independent matrices. Still, $t-$independence is a quite inflexible structural requirement, and our question remains of interest for matrices that are not $t-$independent.  We will come back to this matter in Section \ref{sec.results}.

\subsubsection{The FKP perturbation survives in rough domains and in degenerate elliptic theories} In the past few decades there has been a lot of interest in the Dirichlet problem on domains satisfying weak topological and geometric assumptions. A thorough review of this area is outside our scope,  but some highlighted works include \cite{jk1, sem89, dj,  bl04, hm2}. Of course, one immediately wonders whether the FKP perturbation theory holds in these more general domains. Along these lines, Milakis-Pipher-Toro \cite{mpt14} obtained the analogue of the FKP perturbation result for the bounded chord-arc domains (these have quantitative openness both in the interior and exterior of the domain, as well as quantitative path-connectedness, and their boundary is $(n-1)-$Ahlfors-David regular; see Section \ref{S2} for precise statements). They also obtained the stability of the $RH_p$ condition if the measure on $\Omega$ with density $\rho(A,A_0)^2\delta^{-1}$ is a Carleson measure with small enough norm (depending on $p$ and $A_0$) (see also \cite{esc} and \cite{mt10}).

It is known that the bounded chord-arc domains  have uniformly rectifiable boundaries \cite{dj, hmu}.  Cavero-Hofmann-Martell \cite{chm} have proved that the FKP perturbation theory holds also for real symmetric operators on the more general 1-sided chord-arc domains (quantitative openness and quantitative path-connectedness inside the domain and $(n-1)-\ADR$ boundaries) (see Section \ref{S2}). Their method relies on an extrapolation of Carleson measures technique developed by Lewis and Murray \cite{lm}, which was first used to give an alternate proof of the FKP perturbation result by Hofmann and Martell \cite{hm1}. The technique makes heavy use of sawtooth domains and a Dahlberg-Jerison-Kenig projection lemma which allows one to compare measures on the sawtooth domain to their projections on the original boundary. A year later, Cavero-Hofmann-Martell-Toro \cite{chmt} devised a different method of proof for the FKP perturbation result in the same setting of 1-sided CAD (and extending to the non-symmetric case), using a generalization of a result of Kenig-Kirchheim-Pipher-Toro \cite{kkipt} that weak-$BMO$ solvability of $L$ implies the $A_{\infty}$ property for the elliptic measure of $L$ (in fact, this is a characterization \cite{chmt, hmtbook}).

The state-of-the-art for the elliptic operators satisfying (\ref{eq.elliptic}) lies in the  article of  Akman-Hofmann-Martell-Toro \cite{ahmt}, where they  generalize  the FKP perturbation theory to the situation of uniform domains satisfying a capacity density condition. Since the $(n-1)-$dimensional Hausdorff measure of the boundary of the domain need not be $\ADR$ (indeed, it could potentially be locally infinite), their perturbation result is stated among the elliptic measures only, with no reference to an underlying surface measure. Thus, a main result of theirs reads as follows: Suppose that  $\Omega\subset\bb R^n$ is a bounded (for simplicity, but they consider unbounded domains too) uniform domain satisfying the CDC and fix $X_0$  in the ``center'' of $\Omega$ (for instance, $X_0$ can be any Corkscrew point of a ball with radius $\diam(\partial\Omega)/2$). Let $L_0$, $L$ be two elliptic operators with associated matrices $A_0$, $A$, associated elliptic measures $\omega_0$, $\omega$, and associated Green's functions $G_0$, $G$, respectively. If
\begin{multline}\label{eq.ahmt}
\sup_{r\in(0,\diam(\partial\Omega))}\sup_{x\in\partial\Omega}g(x,r)\\=\sup_{r\in(0,\diam(\partial\Omega))}\sup_{x\in\partial\Omega}\frac1{\omega^{X_0}_0(B(x,r)\cap\partial\Omega)}\dint_{B(x,r)\cap\Omega}\rho(A,A_0)(Y)^2\frac{G_0(X_0,Y)}{\delta(Y)^2}\,dY<+\infty,
\end{multline} 
then $\omega\in A_{\infty}(\omega_0)$ (see Definition \ref{def.ainfty}).   They consider the expression $g(x,r)$ based on an analogous one used as an intermediate step in \cite{fkp}. Furthermore, it is shown that  if there exists a doubling measure $\sigma$  on $\partial\Omega$ such that $\omega_0\in A_{\infty}(\sigma)$ and $\rho(A,A_0)$ satisfies (\ref{eq.ahmt}), then $\rho(A,A_0)$ also satisfies (\ref{eq.hcm}) (with $h(x,r)$ defined using $\sigma$), and $\omega\in A_{\infty}(\sigma)$. In this way, we see that the results of \cite{ahmt} do properly generalize the perturbation results of \cite{fkp}. They are also able to generalize the area integral condition of Fefferman \cite{fe89}, so that if $\n A(\rho(A,A_0))\in L^{\infty}(\partial\Omega,\omega_0)$, then $\omega\in A_{\infty}(\omega_0)$.  They also show that  ``small constant'' assumptions lead to $\omega\ll\omega_0$ and $\frac{d\omega}{d\omega_0}\in RH_p(\omega_0)$. 

More recently, the second author of this paper, together with Svitlana Mayboroda   \cite{mp20}, has used the extrapolation of Carleson measures technique to obtain an analogue of the Carleson perturbation result for the degenerate elliptic operators of David, Mayboroda, and the first author of this paper \cite{dfm1}. Adequate ``small constant'' analogues that preserve the $RH_p$ property are also obtained. A main difficulty in this work is the proper maneuvering of  mixed-dimensional sawtooth domains, for which an elliptic PDE theory is verified using the  axiomatic methods in \cite{dfm20}. The \cite{dfm1} operators are adapted to study domains with low-dimensional $\ADR$ boundaries; for a brief overview on the available results in this direction, see the introduction in \cite{mp20}.

Given the robustness of the FKP perturbations on rough domains and the degenerate elliptic theory, it is natural to ask if the theory holds even in the axiomatic setting of \cite{dfm20}, which allows for domains with mixed-dimensional boundaries, or domains with boundary measures given by non-trivial weights, or even atoms (see more in Section \ref{sec.examples}). 

\begin{question}\label{q4} Does the FKP perturbation theory hold in the axiomatic setting of \cite{dfm20}?
\end{question}

\subsection{Main results}\label{sec.results} We now discuss our contributions, which will give partial or full answers to the questions posed above. 

The domains which we consider are described fully in Section \ref{subpde}, but here let us give a quick review. We assume that our domains $\Omega\subset\bb R^n, n\geq2$ are 1-sided NTA domains (that is, they have the interior Corkscrew and the Harnack Chain properties, see Definition \ref{def1NTA}), and that they are paired with a positive doubling measure $dm=w\,dX$ on $\Omega$ such that $w\in L^1_{\loc}(\Omega,dX)$ and such that $(\Omega,m)$ has an $L^2-$Poincar\'e inequality on interior balls. The weight $w$ is tailored to the study of the boundary of $\Omega$.  We ask our operators $L = - \div A \nabla$ to satisfy an elliptic and boundedness condition that matches the behavior of $w$, that is, for almost every $X\in \Omega$, we require the existence of $C_L >0$ such that
\begin{equation} \label{defellipticw}
A(X) \xi \cdot \xi \geq (C_L)^{-1} w(X) |\xi|^2 \quad \text{ for } \xi \in \R^n,
\end{equation}
and
\begin{equation} \label{defboundedw}
|A(X) \xi \cdot \zeta| \leq C_L w(X) |\xi||\zeta| \quad \text{ for } \xi,\zeta \in \R^n.
\end{equation}
If we write $L$ as $-\div( w\A\nabla)$, we can remove the dependence on $w(X)$ in \eqref{defelliptic}--\eqref{defbounded} and recover the classical elliptic and boundedness conditions
\begin{equation} \label{defelliptic}
\A(X) \xi \cdot \xi \geq (C_L)^{-1} |\xi|^2 \quad \text{ for } \xi \in \R^n
\end{equation}
and
\begin{equation} \label{defbounded}
|\A(X) \xi \cdot \zeta| \leq C_L |\xi||\zeta| \quad \text{ for } \xi,\zeta \in \R^n.
\end{equation}
One may prefer to pick a   second order  operator $L$ first, and then think of $m$ as a way to describe the degeneracies of $L$.   In particular, the case where $L$ is uniformly elliptic and bounded in the classical sense means that $m$ is the Lebesgue measure on $\Omega$, and vice versa.  Finally, we  assume that on our domains $(\Omega,m)$ there is a robust theory for the ``elliptic'' operators in the sense of (\ref{defellipticw}-\ref{defboundedw}). The details of the theory that we assume are laid out in Definition \ref{def.pdef}; in summary, we require boundary H\"older continuity of solutions, the existence and uniqueness of doubling elliptic measures giving the appropriate representation formula for solutions to the continuous Dirichlet problem, and a weakly-defined Green's function.

The domains $(\Omega,m)$ described above are denoted as PDE-friendly domains. We mention several examples in Section \ref{sec.examples}, but a few of them are the 1-sided NTA domains satisfying the capacity density condition, the low-dimensional Ahlfors-David regular domains of \cite{dfm1}, and mixed-dimensional sawtooth domains as in \cite{mp20}.

In the rest of the article, for any $X\in \Omega$, $\delta(X)$ is given as in (\ref{eq.delta}) and $B_X$ denotes $B(X,\delta(X)/4)$. When $x\in \partial \Omega$ and $r>0$, we write $B(x,r)$ for the open ball in $\R^n$ and $\Delta(x,r)$ for the boundary ball $B(x,r) \cap \partial \Omega$. 
Note that the radius $r$ and center $x$ of a boundary ball are not necessary unique, but by a slight abuse of notation, a boundary ball $\Delta$ will mean either a triple $(x,r,\Delta(x,r))$, or just the set $\Delta(x,r)$.  
The truncated area integral $\n A$ and the non-tangential maximal function $N$ are constructed as follows: 
\begin{equation} \label{defA}
\n A^r(f)(x) = \Big(\dint_{\gamma^r(x)} |f(X)|^2 \frac{dm(X)}{m(B_X)}\Big)^\frac12 \qquad \text{ for $r>0$, $x\in \partial \Omega$, and $f\in L^2_{\loc}(\Omega,m)$}
\end{equation}
and
\begin{equation} \label{defN}\nonumber
N^r(f)(x) = \sup_{\gamma^r(x)} |f| \qquad \text{ for $r>0$, $x\in \partial \Omega$, and $f\in C(\Omega)$,}
\end{equation}
where
\begin{equation} \label{defgamma}\nonumber
\gamma^r(x) = \{X \in \Omega, \, |X-x| \leq 2\delta(X) \leq 2r\}.
\end{equation}

\begin{definition}[Doubling family of measures]
We say that a family $\omega = \{\omega^X\}_{X\in \Omega}$ of Borel measures is \emph{doubling} if there exists a constant $C>0$ such that for $x\in  \partial \Omega$, $r>0$, and $X\in \Omega \setminus B(x,4r)$, we have 
\begin{equation} \label{defdoublingf}
\omega^X(B(x,2r)) \leq C \omega^X(B(x,r)).
\end{equation}
The measure $\sigma$ is doubling if \eqref{defdoublingf} is verified with $\sigma$ instead of $\omega^X$. 
\end{definition}

Below, we give a meaning to saying that an $L^2_{\loc}(\Omega,m)$ function satisfies a Carleson measure property.

\begin{definition}[Carleson measure condition] If $\omega = \{\omega^X\}_{X\in \Omega}$ is a doubling family of measures on $\partial \Omega$, we say that a function $f\in L^2_{\loc}(\Omega,m)$ satisfies the \emph{$\omega$-Carleson measure condition} if there exists $M>0$ such that for any $x\in \partial \Omega$,  any  $ r>0$,  and any $X\in \Omega\backslash B(x,2r)$\footnote{Note that if $\diam \Omega < +\infty$ and $r$ is large, then $\Omega\backslash B(x,2r)=\emptyset$ and hence \eqref{defCM} is automatically true (by convention). If $\diam \partial \Omega < +\infty$ and $r$ large, then $\Omega\backslash B(x,2r)\neq\emptyset$ and this definition also makes sense; if we want to weaken the definition to $r\in (0,\diam\partial \Omega)$, then we fall in the situation presented in Subsection \ref{SS2.3}.},
 we have 
\begin{equation} \label{defCM}
\int_{\Delta(x,r)} |\n A^r(f)(y)|^2 \, d\omega^X(y) \leq M \omega^X(\Delta(x,r)).
\end{equation}
We write $f\in KCM(\omega)$  to say that $f$ satisfies the $\omega$-Carleson measure condition and $f\in KCM(\omega,M)$ if we want to refer to the constant in \eqref{defCM}.   We will  often  use $M_f$ for the smallest admissible constant in (\ref{defCM}),  and we call it the \emph{Carleson norm} of $f$.  If $\sigma$ is simply a measure on $\partial\Omega$, the $\sigma$-Carleson measure condition means \eqref{defCM} where $\omega^X$ is replaced by $\sigma$. At last, we shall also need the following variant of the Carleson measure condition. We say that  
\begin{equation} \label{defwtCM}\nonumber
f\in {KCM}_{\sup}(\omega,M) \quad \text{ if$_{def}$ } \quad X \mapsto \sup_{B_X} |f| \in KCM(\omega,M).
\end{equation}
\end{definition}

We note that, in settings where the Green  function $G(X,Y)$ is defined and can be compared to the elliptic measure (via a suitable estimate of the form (\ref{tcp18})), our condition (\ref{defCM}) is equivalent to the finiteness of the expression
\begin{equation}\label{eq.green}
\sup_{r>0}\sup_{x\in\partial\Omega}\sup_{X\in\Omega\backslash B(x,2r)}\frac1{\omega^X(\Delta(x,r))}\dint_{B(x,r)\cap\Omega}|f(Y)|^2\frac{G(X,Y)}{\delta(Y)^2}\,dm(Y)
\end{equation}
via Fubini's theorem. By comparing  (\ref{eq.green}) to (\ref{eq.ahmt}), we see that our Carleson measure condition is a  reformulation of analogue Carleson measure conditions considered in \cite{fkp} and \cite{ahmt}. Moreover, since $\n A^r\leq\n A$ (where $\n A$ is the area integral with no truncation), our condition (\ref{defCM}) readily captures the same results under (an analogue of) the stronger $L^{\infty}$  assumption on the area integral (\ref{eq.areaint}); this last observation had essentially been made already in \cite[Chapter 3]{ahmt}.
 
Now, we define $A_{\infty}-$absolute continuity  among doubling families of  measures.

\begin{definition}[$A_{\infty}$ for families of measures]\label{def.ainfty} If $\omega_0 = \{\omega^X_0\}_{X\in \Omega}$ and $\omega_1 = \{\omega^X_1 \}_{X\in \Omega}$ are two doubling families of measures on $\partial \Omega$, then we say that $\omega_1$ is \emph{ $A_\infty$-absolutely continuous with respect to $\omega_0$} - or $\omega_1 \in A_\infty(\omega_0)$ for short - if, for any $\xi>0$, there exists $\zeta>0$ such that for any boundary ball $\Delta := \Delta(x,r)$, any $X\in \Omega \backslash B(x,2r)$, and any Borel set $E \subset \Delta$, we have that
\begin{equation} \label{defAinfty}
    \frac{\omega_1^X(E)}{\omega_1^X(\Delta)} < \zeta\quad \text{ implies }\quad \frac{\omega_0^X(E)}{\omega_0^X(\Delta)} < \xi. 
\end{equation}
If $\sigma_0$ or $\sigma_1$  are  measures, then we replace $\omega^X_0$ by $\sigma_0$ or/and $\omega^X_1$ by $\sigma_1$ in \eqref{defAinfty}. 
\end{definition}

Our first main theorem links a bound on the oscillations of bounded solutions to $A_\infty$. The result is the analogue in our setting of \cite[Theorem 1.1 (a)$\implies$(b)]{chmt} or \cite[Theorem 4.1]{kkipt}.

\begin{theorem}[Weak-$BMO$ solvability implies $A_{\infty}$]\label{Main1}
Let $(\Omega,m)$ be a PDE friendly domain (see Definition \ref{def.pdef}). Let $L= - \div A\nabla$ be an elliptic operator satisfying \eqref{defellipticw} and \eqref{defboundedw}, and construct the elliptic measure $\omega:= \{\omega^X\}_{X\in \Omega}$ as in \eqref{defhm}.  Let $\sigma$ be a doubling measure or doubling family of measures on $\partial \Omega$. 

If there exists $M>0$ such that,  for any Borel $E \subset \partial \Omega$,  the solution $u_E$ constructed as $u_E(X) := \omega^X(E)$ satisfies
\begin{equation} \label{Main1a}
\delta \nabla u_E \in KCM(\sigma,M),
\end{equation}
then $\omega \in A_\infty(\sigma)$.
\end{theorem}

In fact, we prove stronger local analogues; see Lemma \ref{Main1loc} and Corollary \ref{Main11}.

Our second main theorem states that Carleson perturbations of an elliptic operator perserve the $A_\infty$-absolute continuity, via an $S<N$ estimate.   However, we give a much broader sense to Carleson perturbations than what was found previously in the literature, and that will be our contribution to the answer of Question \ref{q1}.

\begin{definition}[Generalized Carleson perturbations] \label{defGCP}
Let $L_0 = -\div (w\A_0 \nabla)$ and $L_1=-\div (w\A_1 \nabla)$ be two operators satisfying \eqref{defelliptic}--\eqref{defbounded}, and let $\omega_0 = \{\omega^X_0\}_{X\in \Omega}$ be the elliptic measure of $L_0$ constructed in \eqref{defhm}.

\noindent We say that $L_1$ is an {\bf additive Carleson perturbation} of $L_0$ if 
\[ |\A_1 - \A_0| \in KCM_{\sup}(\omega_0).\]
We say that $L_1$ is a {\bf scalar-multiplicative Carleson perturbation} of $L_0$ if there exists a scalar function $b$ such that $C^{-1} \leq b \leq C$ for some $C>0$ and 
\[ \A_1 = b\A_0, \qquad \text{ and } \delta|\nabla b| \in KCM(\omega_0).\]
The operator $L_1$ is an {\bf antisymmetric Carleson perturbation} of $L_0$ if there exists a  bounded, antisymmetric matrix-valued function $\mathcal T$ such that  
\[ \A_1 = \A_0 + \mathcal T ,\qquad \text{ and } \delta w^{-1} |\div w\mathcal T| \in KCM(\omega_0)\]
where $\div(\T)$ is the vector obtained by taking the divergence of each column of $\T$. At last, $L_1$ is a {\bf (generalized) Carleson perturbation} of $L_0$ if there exists a matrix-valued function $\mathcal C$, a scalar function $b$, and an antisymmetric matrix-valued function $\mathcal T$ such that
\[ 
|\mathcal C| \in  KCM_{\sup}(\omega_0), \quad \delta| \nabla b| + \delta w^{-1}|\div(w \T)| \in KCM(\omega_0), \quad \text{and}  \quad \A_1 = b(\A_0 + \mathcal C + \T),
\]
and the {\bf norm of the Carleson perturbation} is the smallest value $K>0$ such that $ |\mathcal C| \in KCM_{\sup}(\omega_0,K)$ and  $\delta|\nabla b|/b + \delta w^{-1} |\div (w\T)| \in KCM(\omega_0,K)$.  
\end{definition}

Note that the additive Carleson perturbation is what was known as the Carleson perturbation in earlier articles, and so we extended the notion of Carleson perturbation to the `scalar-multiplicative' and `antisymmetric' perturbations. These  last two types of perturbation  can be seen (at least formally) as {\bf drift Carleson perturbations} via the following well-known transformations:
\begin{equation}\label{eq.scalarmult}
L_1 := -\div (w b\A_0 \nabla) = - b \div(w \A_0 \nabla) - w (\A_0)^T \nabla b \cdot \nabla = bL_0 - w (\A_0)^T \nabla b \cdot \nabla
\end{equation}
and
\begin{equation}\label{eq.antisym}
L_1 := -\div (w[\A_0 + \T] \nabla)  = L_0 - \div( w\T) \cdot \nabla.
\end{equation}
On the other hand, note that our perspective allows us to consider these perturbations without a priori constructing an elliptic theory for operators with drift terms. The scalar-multiplicative and antisymmetric perturbations  are interesting because they are perturbations that can significantly change the coefficients of the initial matrix $\A_0$ in a neighborhood of the boundary, thus answering Question \ref{q3}.
They also appeared naturally in previous works. In \cite{dm} and \cite{fen1}, the authors proved that, when the boundary is a uniformly rectifiable set of dimension $d<n-1$, the elliptic measure associated to the operators $L_\beta = -\div [D_\beta]^{d+1-n} \nabla$ is $A_\infty$-absolutely continuous with respect to the $d$-dimensional approach (see \cite{dm}, \cite{fen1} for the definitions of uniformly rectifiable and $D_\beta$); the proof in \cite{fen1} relies on the fact $L_\beta$ are scalar-multiplicative Carleson perturbations of each other.
Theorem 1.6 in \cite{chmt} states a particular case of the following assertion, which is an easy consequence of our Theorem \ref{Main2} below: if $L^*$ is an antisymmetric Carleson perturbation of $L$, then $\omega_{L^*} \in A_\infty(\omega_L)$, and the elliptic measure of the self-adjoint operator $L_{s} =(L+L^*)/2$ belongs to the same $A_\infty$ class than $\omega_{L}$ and $\omega_{L^*}$. The idea of taking Carleson perturbations in the drift term has also appeared before \cite{hl, kp3}, but to the best of our knowledge, the present article is the first time that drift Carleson perturbation are used to extend the class of transformations of the elliptic matrix $A$ that perserves the $A_\infty$-absolute continuity.

\begin{theorem}[$S<N$ is preserved by Carleson perturbations]\label{Main2}
Let $(\Omega,m)$ be a PDE friendly domain (see Definition \ref{def.pdef}). Let $L_0= - \div w\A_0\nabla$ and $L_1 = - \div w\A_1 \nabla$ be two elliptic operators satisfying \eqref{defelliptic} and \eqref{defbounded}, and construct the elliptic measures $\omega_0:= \{\omega_0^X\}_{X\in \Omega}$ and $\omega_1:= \{\omega_1^X\}_{X\in \Omega}$  as in \eqref{defhm}. 

If $L_1$ is a (generalized) Carleson perturbation of $L_0$, then for any $x\in \partial \Omega$, any $r\in (0,\diam \Omega)$, any Corkscrew point $X$ associated to $(x,r)$, and any weak solution $u$ to $L_1 u=0$, we have that
\begin{equation} \label{Main2b}
\int_{\Delta(x,r)} |\n A^r(\delta \nabla u)|^2 d\omega_0^X \leq C \int_{\Delta(x,2r)} |N^{2r}(u)|^2 d\omega^X_0,
\end{equation}
with a constant $C>0$ that depends only on the dimension $n$, $C_{L_0}$, $C_{L_1}$, the norm of the Carleson perturbation, and the constants in the PDE friendly properties of $(\Omega,m)$.

In particular, \eqref{Main1a} holds with $\sigma = \omega_0$, and hence $\omega_1 \in A_\infty(\omega_0)$.
\end{theorem}

Via different methods, a local $S<N$ result (which works even in more general $L^q$ settings) has been obtained in \cite{ahmt} for the 1-sided NTA domains satisfying the capacity density condition.   We could also obtain the same result from \cite{ahmt} by applying a good-$\lambda$ argument to \eqref{Main2b}, but  we do not need those bounds for the present paper and decided to postpone them for a future article.

It is well know that    $A_\infty$  is an equivalence relationship (see     \cite{gr}), which means that Theorem \ref{Main2} would also hold if we assume that $L_0$ is a Carleson perturbation of $L_1$ (which is {\em a priori} different from saying that $L_1$ is a Carleson perturbation of $L_0$, since the Carleson measure condition depends on the operator before perturbation). However,  by combining Theorem \ref{Main2} with the theorem below, we obtain that our notion of `Carleson perturbations of elliptic operators' is actually an equivalence relationship, as expected. 

\begin{theorem}[$A_{\infty}$ implies transitivity of $CM$]\label{Main3} Let $(\Omega,m)$ be a PDE friendly domain (see Definition \ref{def.pdef}), and for $i\in \{0,1\}$, let $\mu_i$ be either an elliptic measure [$\mu_i = \{\omega_i^X\}_{X\in\Omega}$] or a doubling measure [$\mu_i = \sigma_i$]  on $\partial\Omega$. If $\mu_1\in A_{\infty}(\mu_0)$, then
\begin{equation}\label{eq.transcm}
f \in KCM(\mu_0) \quad\text{if and only if }\quad f \in KCM(\mu_1),\qquad\text{for each }f\in L^2_{\loc}(\Omega,m).
\end{equation} 
\end{theorem} 

For a local analogue of the above result, see Lemma \ref{lm.jn2}.  We actually can prove a characterization of $A_{\infty}$ via the property (\ref{eq.transcm}), see Corollary \ref{cor.char} below. Theorem \ref{Main3} can be seen as analogue of the John-Nirenberg lemma (which is for $BMO$ functions) adapted to Carleson measures and $A_\infty$ weights.  The result is an extension to our setting of  \cite[Lemma 3.8]{hmm20}, which itself is a modification of  John-Nirenberg type inequalities proved in \cite[Lemma 10.1]{hmay}, \cite[Lemma 2.14]{ahlt}, and \cite[Lemma A.1]{mmm}, although our method of proof is different.
Since the condition ${KCM}_{\sup}(\omega_i)$ is only $KCM(\omega_i)$ applied to a transformation of $f$, we have in particular that if $\omega_1\in A_{\infty}(\omega_0)$, then $f \in {KCM}_{\sup}(\omega_0) \Leftrightarrow f \in {KCM}_{\sup}(\omega_1)$. Lastly, see Lemma \ref{lm.jn2} for a local version.

Let us emphasize that none of our proofs rely on the construction of sawtooth domains on PDE friendly domains, nor do they rely on the extrapolation theory of Carleson measures. Indeed, it is not clear to us that sawtooth domains of PDE friendly domains are themselves PDE friendly. Even if they were, the   construction of, and verification of   PDE friendly axioms on sawtooth domains of some rough domains are long and difficult tasks \cite{hm2, mp20}. 
Our method resembles loosely that of the recent paper \cite{chmt}, where an analogue of Theorem \ref{Main1} is used to extend the FKP (additive) perturbation theory to the case of 1-sided chord-arc domains, but they also use sawtooth domains.

The rest of the article will be divided as follows.   In the rest of the introduction, we give some applications of our three theorems (Theorem \ref{Main1}, Theorem \ref{Main2}, Theorem \ref{Main3}).  In Section \ref{S2}, we present the assumptions for the PDE friendly domains and examples of these domains. In Section \ref{S3}, we recall the theory of $A_\infty$-weights that we need for our proof, and moreover, we prove Theorem \ref{Main3}. Section \ref{S4} and Section \ref{S5} are devoted to the proofs of  Theorem \ref{Main1} and Theorem \ref{Main2}, respectively.

\subsection{Applications  of  main results}\label{sec.cor} Let us present several implications of our theorems.

First, a straightforward consequence of Theorems \ref{Main1}, \ref{Main2}, and \ref{Main3} is the fact that if the elliptic measure $\omega_0$ is already $A_\infty$-absolutely continuous with respect to a doubling measure $\sigma$, and $L_1$ is a Carleson perturbation of $L_0$, then the $A_\infty(\sigma)$ absolute continuity is transmitted to $\omega_1$. Thus, our results extend the FKP perturbation theory to PDE friendly domains, giving in particular  a positive answer to our Question \ref{q4} (since the axiomatic setting of \cite{dfm20} satisfies the assumptions of the PDE friendly domains).

\begin{corollary}[An extension of the FKP perturbation result to PDEF domains]\label{Main4} Let $(\Omega,m)$ be a PDE friendly domain (see Definition \ref{def.pdef}), and let $\sigma$ be a doubling measure on $\partial \Omega$. Consider an elliptic operator  $L_0:= -\div (w \mathcal  A_0 \nabla)$ such that $\omega_0 \in A_\infty(\sigma)$. Assume that the elliptic operator $L_1:= -\div (w \mathcal  A_1 \nabla)$ is a Carleson perturbation of $L_0$ in the sense that there exist  a matrix $\mathcal C$, a scalar $b$, and an antisymmetric matrix $\mathcal T$ such that $\mathcal A_1 = b(\mathcal A_0 + \mathcal C + \mathcal T)$, and
\[ 
|\mathcal C| \in  KCM_{\sup}(\sigma), \quad \text{ and } \quad  \delta| \nabla b| + \delta w^{-1}|\div(w \T)| \in KCM(\sigma).
\]
Then $\omega_1 \in A_\infty(\sigma)$.
\end{corollary}

\bp
Since $\omega_0 \in A_\infty(\sigma)$, Theorem \ref{Main3} gives that $L_1$ is a generalized Carleson perturbation of $L_0$, as given in Definition \ref{defGCP}, hence as needed for Theorem \ref{Main2}. Applying Theorem \ref{Main2} and then Theorem \ref{Main1} yields the desired $\omega_1 \in A_\infty(\sigma)$.
\ep

Moreover, our theory gives

\begin{corollary}[Equivalence of $A_{\infty}$ and weak-$BMO-$solvability] \label{Main5} Let $(\Omega, m)$ be a PDE friendly domain (see Definition \ref{def.pdef}). Let $L$ and $\omega$ as in Theorem \ref{Main1}, and take a doubling measure $\sigma$ on $\partial\Omega$. The following are equivalent:
\begin{enumerate}[(i)]
\item\label{bmo1} $\omega \in A_\infty(\sigma)$.
\item\label{bmo2} the Dirichlet problem to $L u = 0$ is weak-$BMO(\sigma)$ solvable; that is, there exists $M>0$ such that, for any Borel $E \subset \partial \Omega$, the solution $u_E$ constructed as $u_E(X) := \omega^X(E)$ satisfies
\begin{equation} \label{Main5a}\nonumber
\delta \nabla u_E \in KCM(\sigma,M).
\end{equation}
\end{enumerate}
\end{corollary}

\bp The implication \ref{bmo2} $ \Rightarrow$ \ref{bmo1} is a consequence of Theorem \ref{Main1}. Since $L$ is an $\omega$-Carleson perturbation of itself, \ref{bmo1} $ \Rightarrow$ \ref{bmo2} follows from Theorems \ref{Main2} and \ref{Main3}. 
\ep

After   the first version of our paper was posted online, we learned that  Cao, Dom\'inguez, Martell, and Tradacete were about to finish an article (see \cite{cdmt}) with a result similar to our Corollary \ref{Main5}, 
using roughly the same techniques as the ones we used. 
They worked on the specific case of 1-sided NTA domains satisfying the capacity-density condition, but they gave many more characterizations of   the $A_\infty$ property of the elliptic measure than us. 
In particular, they prove that the full $BMO$ solvability, and  the $S<N$ estimate in some $L^q$ (with $S$ the conical square function), are equivalent to the $A_{\infty}$ property among elliptic measures. 
Note however that the main interest of \cite{cdmt} differs from ours, in that they focused on criterions of $A_\infty$, while we were primarily interested in extending the notion of Carleson perturbations that preserves $A_\infty$, an issue which is not considered in \cite{cdmt}. 

Next, we show that our Theorems \ref{Main1} and \ref{Main2} yield a certain converse to Theorem \ref{Main3}, which gives a new characterization of $A_{\infty}$ among elliptic measures, via the transitivity of the Carleson measure property. Thus, the following corollary gives an answer to Question \ref{q2}, regarding the connection between $A_{\infty}$ and Carleson measures. This characterization of $A_{\infty}$ seems new to us, and it does not appear in \cite{cdmt} either.  

\begin{corollary}[$A_{\infty}$ is equivalent to transitivity of $CM$, for elliptic measures]\label{cor.char} Let $(\Omega,m)$ be a PDE friendly domain (see Definition \ref{def.pdef}), and let $\{\omega_0^X\}_{X\in\Omega}$, $\{\omega_1^X\}_{X\in\Omega}$ be two elliptic measures on $\partial\Omega$. Then, $\omega_1\in A_{\infty}(\omega_0)$ if  and only if 
	\begin{equation}\label{eq.transcm1}
	f \in KCM(\omega_1) \quad\text{implies that }f \in KCM(\omega_0),\qquad\text{for each }f\in L^2_{\loc}(\Omega,m).
	\end{equation} 
\end{corollary} 

\noindent\emph{Proof.} The ``only if'' direction is immediate from Theorem \ref{Main3}. Now suppose that  (\ref{eq.transcm1}) holds.   Let $E\subset\partial\Omega$ be an arbitrary Borel set, and write $u_1(X):=\omega_1^X(E)$. Note that $\delta|\nabla u_1|\in L^2_{\loc}(\Omega, m)$. According to Theorem \ref{Main2}, we have (\ref{Main2b}), which implies in particular that $\delta\nabla u_1\in KCM(\omega_1)$. By   hypothesis, it follows that $\delta\nabla u_1\in KCM(\omega_0)$. Since $E$ was arbitrary, then Theorem \ref{Main1} allows us to conclude that $\omega_1\in A_{\infty}(\omega_0)$.\hfill{$\square$}

It seems to us that  Corollary \ref{cor.char} has not been known even in the classical settings of the half-space or the unit ball. It is not clear that the FKP characterization of (classical) $A_{\infty}$ via a Carleson measure condition \eqref{eq.fkpchar} immediately implies a suitable analogue of our Corollary \ref{cor.char}. On the other hand, we emphasize that our characterization is proved only among elliptic measures; whether Corollary \ref{cor.char} holds for general doubling measures is an open question, even in the case of the half-space.
 
We also remark that the scalar-multiplicative Carleson perturbations contain the scalar subclass of Dahlberg-Kenig-Pipher operators \ref{item.linfty}-\ref{item.cm}. More precisely, it is easy to see that if $A=bI$ is a matrix satisfying the ellipticity and boundedness conditions (\ref{eq.elliptic}) and the DKP conditions \ref{item.linfty}-\ref{item.cm}, then $b$ verifies the assumptions
\begin{equation}\label{eq.bass}
C^{-1}\leq b\leq C,\qquad\text{and}\qquad \delta\nabla b\in KCM(\sigma),
\end{equation}
where $\sigma$ is the surface measure; and on the other hand, if $b$ verifies (\ref{eq.bass}),  then $A=bI$ satisfies \ref{item.cm}. By seeing this subclass as a scalar-multiplicative perturbation from the Laplacian $-\Delta$, we are able to obtain, for instance, alternate proofs of difficult results for the scalar subclass of DKP operators, which have recently been shown for the full generality of DKP operators. As a matter of fact, our result for the scalar operators goes slightly beyond that of the DKP operators, as we do not have to assume the boundedness condition on the gradient \ref{item.linfty}. Pointedly, consider

\begin{corollary}[A free boundary result for scalar DKP operators]\label{cor.fb} Let $\Omega\subset\bb R^n$, $n\geq3$, be a uniform (that is, 1-sided NTA) domain with $(n-1)-$Ahlfors-David regular boundary (see Definition \ref{def.adr}), and set $\sigma=\m H^{n-1}|_{\partial\Omega}$. Let $b$ be a function on $\Omega$ verifying $C^{-1}\leq b\leq C$ and $\delta\nabla b\in KCM(\sigma)$. Then the following are equivalent.
\begin{enumerate}[(i)]
	\item The elliptic measure $\omega_L$ associated with the operator $L=-\dv b\nabla$ is $A_{\infty}$ with respect to $\sigma$.
	\item $\partial\Omega$ is uniformly rectifiable.
	\item $\Omega$ is a chord-arc domain.
\end{enumerate}
\end{corollary}

\noindent\emph{Sketch of proof.} For definitions of uniform rectifiability and chord-arc domain, see for instance \cite{hmmtz2}. For $L\equiv-\Delta$, then the above equivalences are known \cite{ahmnt, ahmmt}; in particular, (ii) and (iii) are equivalent, and either imply (i) with $L=-\Delta$.

We show (i)$\implies$(ii); the converse has a similar proof. Say that $L=-\dv b\nabla$ and $b$ has the described properties, and suppose that $\omega_L\in A_{\infty}(\sigma)$. Then $b^{-1}$ also has the same properties as $b$; that is, $C^{-1}\leq\frac1b\leq C$ and $|\delta\nabla(1/b)|=\delta\nabla(b)/b^2\in KCM(\sigma)$. Since $-\Delta=-\dv(b^{-1}b\nabla)$, then $-\Delta$ is a scalar-multiplicative Carleson perturbation of $L$, and by Corollary \ref{Main4} it follows that $\omega_{-\Delta}\in A_{\infty}(\sigma)$.  Thus we have that $\partial\Omega$ is uniformly rectifiable.\hfill{$\square$}

It is clear that Corollary \ref{cor.fb} also holds with $L$ being any generalized Carleson perturbation from $-\Delta$. The above result for functions $b$ which also verify that $\delta\nabla b\in L^{\infty}(\Omega)$ is a particular case of the recent free boundary result for the DKP operators shown by Hofmann-Martell-Mayboroda-Toro-Zhao \cite{hmmtz}; our method of Carleson perturbations allows us to dispense with the aforementioned boundedness condition.

Finally, let us consider an application of our theory to the study of elliptic measures on purely unrectifiable sets\footnote{The authors would like to thank Max Engelstein for pointing this application out.}. Let $K$ be the Garnett-Ivanov Cantor set \cite{garnettbook, ivanov}, also known as the $4$-corner Cantor set, as defined in Section 3 of \cite{dm2}. The set $K$ is $1-$Ahlfors-David regular with surface measure $\sigma=\n H^1$, and it is known to be purely unrectifiable, and such that the harmonic measure $\omega_{-\Delta}$ and surface measure $\sigma=\n H^1$ are mutually singular.  In \cite{dm2}, G. David and S. Mayboroda constructed an elliptic operator $L_a=-\dv a\nabla$ on $\Omega:=B\backslash K\subset\bb R^2$ such that $\omega_{L_a}\in A_{\infty}(\sigma)$ on $\partial\Omega$, where $B$ is the unit ball in the plane centered at the origin (and $K\subset B$), and $a$ is a certain scalar real-valued function with $1/C\leq a\leq C$ in $\Omega$. It is  not hard to show that the domain $\Omega$ is a 1-sided NTA domain (that is, it has interior Corkscrews and Harnack Chains). Hence $(\Omega,\n L^2)$ is a PDE friendly domain, whence our perturbation theory applies.  We can use our theory to study certain a priori properties of the function $a$ constructed in \cite{dm2}.
	
First, any elliptic operator $L=-\dv A\nabla$ which is a generalized Carleson perturbation of $L_a$ verifies that $\omega_L\in A_{\infty}(\sigma)$ by Corollary \ref{Main4}\footnote{In the case of the classical \emph{additive} Carleson perturbations, this consequence already follows from the perturbation theory in \cite{ahmt}.}. On the other hand, since $\omega_{-\Delta}\perp\sigma$, it follows that the operator $-\Delta$ on $\Omega$ cannot be a generalized Carleson perturbation of $L_a$; in particular, it cannot be a scalar-multiplicative Carleson perturbation. When we put this realization together with Corollary \ref{Main5} and Fubini's theorem, we conclude that 
\[
\sup_{x\in K}\sup_{r\in(0,2)}\frac1{\sigma(\Delta(x,r))}\dint_{B(x,r)}\delta(Y)|\nabla a(Y)|^2\,dY=+\infty.
\]
 In other words, the measure with density $\delta|\nabla a|^2$ is not a Carleson measure.

\section{Hypotheses and elliptic theory} \label{S2}

 Throughout, our ambient space is $\bb R^n$, $n\geq 2$.  We will often write $a\lesssim b$ to mean that there exists a constant $C\geq1$ such that $a\leq Cb$, where $C$ may depend only on certain allowable parameters. Likewise, we write $a\approx b$ if there exists a constant $C\geq1$ such that $\frac1Cb\leq a\leq Cb$. If $(\Omega,\sigma)$ is a measure space and $E\subset\Omega$ is measurable, we will always write $\dashint_E f\,d\sigma=\frac1{\sigma(E)}\int_Ef\,d\sigma$.

\subsection{PDE friendly domains}\label{subpde} In this section we describe the PDE friendly domains and present several examples. First, let us set up some background definitions. Let $n\geq2$ and $\Omega \subset \R^n$ be open.

\begin{definition}[The doubling measure $m$ on the domain] \label{defm} For the remainder of the article, we denote by $m$   a measure on $\Omega$ that satisfies the following  properties:
\begin{enumerate}[(i)]
	\item The measure $m$ is absolutely continuous with respect to the Lebesgue measure; that is, there exists a non-negative weight $w\in L^1_{\loc}(\Omega)$  such that for each Borel set $E\subset\Omega$, $m(E)=\iint_Ew(X)\,dX$.
	\item The measure $m$ is doubling, meaning that there is a constant $C_m\geq1$ such that 
	\begin{equation} \label{mdoubling}
	m(B(X,2r)\cap\Omega) \leq C_{m} m(B(X,r)\cap\Omega) \qquad \text{ for each } X\in \overline{\Omega} \text{ and } r >0.
	\end{equation}
	\item For any open set $D$ compactly contained in $\Omega$,  and any sequence  $\{u_i\}_i\subset C^\infty(\overline {D})$  verifying that  $\iint_D |u_i| \, dm \to 0$ and $\iint_D |\nabla u_i - v|^2\, dm \to 0$ as $i\to \infty$, where $v$ is a vector-valued function in $L^2(D,m)$, we have that $v\equiv 0$.
	\item We assume an $L^2-$Poincar\'e inequality on interior balls:  there exists $C_P$ such that for any ball $B$ satisfying $2B \subset \Omega$ and any function $u\in W^{1,2}(B,m)$, one has
	\begin{equation} \label{defPoincare}
	\fiint_B |u - u_B|^2 \, dm \leq C_P  r \Big( \fiint_B |\nabla u|^2 \, dm \Big)^\frac12,
	\end{equation}  
	where $u_B$ stands for $\dashint_B u\, dm$ and $r$ is the radius of $B$. 
\end{enumerate} 
\end{definition}

Let us briefly discuss our assumptions on $m$. The space $L^2_{\loc}(\Omega,m)$ is not necessarily a space of distributions, meaning that   we may not access the notion of a distributional gradient. However, 
as in \cite{hkmbook} and \cite{dfm20}, the assumption (iii) allows us to construct a notion of gradient $\nabla$ for functions in $L^2_{\loc}(\Omega,m)$, and then we let $W^{1,2}_{\loc}(\Omega,m)$ be the space of functions in $L^2_{\loc}(\Omega,m)$ whose gradient is also in $L^2_{\loc}(\Omega,m)$. It is in this sense that we take the gradient in (\ref{defPoincare}). 

\begin{remark}
As long as the weight $w$ that defines $m$ satisfies the slowly varying property
\begin{equation} \label{mregular}
\sup_{B} w \leq C \inf_B w \qquad \text{ for any ball $B$ such that $2B \subset \Omega$},
\end{equation}
then $L^2(\Omega,m)$ is a space of distributions, and the gradient on $L^2(\Omega,m)$ is the gradient in the sense of distribution. In addition, \eqref{defPoincare} is true. So as long as \eqref{mregular} is verified, we just need to take $m$ such that \eqref{mdoubling} is true.
\end{remark}
 
From there, we can consider the operator $L=-\div (w\A \nabla)$ that satisfies \eqref{defelliptic} and \eqref{defbounded}. We say that $u$ is a \emph{weak solution}  to $Lu = 0$ if $u \in W^{1,2}_{\loc}(\Omega,m)$ and satisfies
\begin{equation} \label{defweak}\nonumber
\dint_\Omega \A \nabla u \cdot \nabla \varphi \, dm = 0 \qquad \text{ for each } \varphi \in C^\infty_c(\Omega).
\end{equation} 
We can deduce the Harnack inequality. 

\begin{lemma}[Harnack inequality, Theorem 11.35 in \cite{dfm20}]  \label{Harnack}
Let $\Omega \subset \R^n$ and $m$ be as in Definition \ref{defm}, and $L=-\div (w\A \nabla)$ satisfy \eqref{defelliptic} and \eqref{defbounded}. If $B$ is a ball such that $2B\subset \Omega$, and if $u\in W^{1,2}_{\loc}(\Omega,m)$ is a non-negative solution to $Lu=0$ in $2B$. Then 
\begin{equation} \label{Harnack1}
\sup_B u \leq C \inf_B u,
\end{equation}
where $C$ depends only on $n$, $C_m$, $C_{P}$, and $C_L$.
\end{lemma} 

Our results are about boundaries, more exactly measures and elliptic measures on the boundary. So, in order to link solutions in $\Omega$ and properties of $\partial \Omega$, we require the domain $\Omega$ to have enough access to the boundary.

\begin{definition}[1-sided NTA] \label{def1NTA}
We say that $(\Omega,m)$ is a  \emph{1-sided} NTA domain  if the following two conditions holds.
\begin{enumerate}
\item[] {\bf Corkscrew point condition} (quantitative openness){\bf.} There exists $c_1\in (0,1)$ such that for any $x\in \partial \Omega$ and any $r\in(0,\diam\Omega)$ we can find $X$ such that $B(X,c_1r) \subset B(x,r) \cap \Omega$.

For $x\in \partial \Omega$ and $r>0$, we say that $X$ is a \emph{Corkscrew point} associated to the couple $(x,r)$  if $
c_1r/100 \leq \delta(X) \leq |X-x| \leq 100r$.
 
\item[] {\bf Harnack chain condition} (quantitative path-connectedness){\bf.} For any $\Lambda \geq 1$, there exists $N_\Lambda$ such that if $X,Y \in \Omega$ satisfy $\delta(X)>r$, $\delta(Y) >r$, and $|X-Y| \leq \Lambda r$, then we can find $N_\Lambda$ balls $B_1, \dots, B_N$ such that $X\in B_1$, $Y \in B_{N_\Lambda}$, $2B_i \subset \Omega$ for $i \in \{1,N_{\Lambda}\}$, and $B_i \cap B_{i+1} \geq 0$ for $i \in \{1,N_{\Lambda}-1\}$.
\end{enumerate}
\end{definition}

\begin{remark} \label{rm.Harnack}
In the Harnack chain condition, we can assume without loss of generality that $X$ is the center of $B_1$, that $Y$ is the center of $B_{N_\Lambda}$, and that $20B_i \subset \Omega$. We may have to increase the value of $N_\Lambda$ but it will still be independent of $X$, $Y$, and $r$.
\end{remark}

At last, for our results to hold, we need   a nice elliptic theory. For the purpose of the article, we shall state the results that we need here, and some geometric settings where they hold.

\begin{definition}[PDE friendly domains]\label{def.pdef}
We say that $(\Omega,m)$ is PDE friendly if $\Omega$ is 1-sided NTA, if $m$ is as in Definition \ref{defm}, and if we have the following elliptic theory. 

Let $L = -\div (w\A \nabla)$ be any second order divergence order operator, where $w$ is the weight in Definition \ref{defm}  and where $A$ is matrix with measurable coefficients which satisfies the ellipticity and boundedness conditions \eqref{defelliptic}--\eqref{defbounded}.

\begin{enumerate}
\item[] {\bf Existence and uniqueness of elliptic measure.} There exists an elliptic measure associated to $L$, which is the only family of probability measures $\omega_L = \{\omega_L^X\}_{X\in \Omega}$ on $\partial \Omega$ such that, for any function $f\in C_c^\infty(\R^n)$, the function $u_f$ constructed as
\begin{equation} \label{defhm}
u_f(X) = \int_{\partial \Omega} f(y) d\omega_L^X(y), \qquad \text{ for } X\in \Omega,
\end{equation}
is continuous on $\overline{\Omega}$, satisfies $u_f = f$ on $\partial \Omega$, and is a weak solution to $Lu=0$.

\item[] {\bf Doubling measure property.} For $x\in \partial \Omega$ and $r>0$, we have that
\begin{equation} \label{dphm1} 
\omega_L^X(\Delta(x,2r)) \leq C \omega_L^X (\Delta(x,r)) 
 \ \ \ \text{ for }
X \in \Omega \setminus 3B,
 \end{equation}
where $\Delta(x,r):=B(x,r)\cap\partial\Omega$, and $C>0$ is independent of $x$, $r$, and $X$, and depends on $L$ only via $C_L$.
 
\item[] {\bf Change of pole.}  Let $x\in \partial \Omega$, $r>0$, and $X$ be a Corkscrew point associated to $(x,r)$. If $E \subset \Delta(x,r)$ is a Borel set, then 
  \begin{equation} \label{CP18} 
C^{-1} \omega_L^{X}(E) \leq  \frac{\omega_L^{Y}(E)}{\omega_L^{Y}(\Delta(x,r))} 
\leq C \omega_L^{X}(E), \ \ \ \text{ for } Y \in \Omega \setminus  B(x,2r),
 \end{equation}
where $C>0$ is independent of $x$, $r$, $X$, $E$ and $Y$,  and depends on $L$ only via $C_L$. 
 
\item[] {\bf H\"older regularity at the boundary.} For any $X\in \Omega$ and any Borel set $E\subset \partial \Omega$, we have
\begin{equation} \label{Holder}
\omega_L^X(E) \leq C \Big( \frac{\delta(X)}{\dist(X,E)} \Big)^\gamma.
\end{equation}
where $C>0$ and $\gamma\in (0,1)$ are independent of $X$ and $E$,  and depend on $L$ only via $C_L$.
 
\item[] {\bf Comparison with the Green function.} Let $X\in \Omega$, write $r$ for $\delta(X)/2$, and take $x \in \partial \Omega$ such that $|X-x| = 2r$. That is, $X$ is a Corkscrew point associated to $(x,2r)$.  There exists a weak solution $G^*_X$ to $L^*u = - \div (w\m A^T \nabla)=0$ in $B(x,r) \cap \Omega$ such that if $y \in \Delta(x,r)$, $s\in (0,r)$, and $Y \in B(x,r) \cap \Omega$ is a Corkscrew point associated to $(y,s)$, we have 
 \begin{equation} \label{tcp18}
C^{-1} \frac{m(B(y,s)\cap \Omega)}{s^2} G^*_X(Y) \leq \omega_L^X(\Delta(y,s)) \leq C \frac{m(B(y,s) \cap \Omega)}{s^2} G^*_X(Y). 
\end{equation}
where $C>0$ is independent of $X$, $y$, $s$ and $Y$, and depends on $L$ only via $C_L$.

Of course, when we write $G^*_X(Y)$, we think of the Green function associated to $L^*$ with pole at $X$. Indeed, in a setting where the   notion of the Green function has been developed, like in \cite{dfm1}, we write $g(X,Y)$ for the Green function associated to $L$ with pole at $Y$, and we set $G^*_X(Y) := g(X,Y)$. In this case, the bounds \eqref{tcp18} are a consequence of \cite[Lemma 15.28]{dfm1} and the fact that $G^*_X$ is a weak solution to $L^*u$ comes from \cite[Lemma 14.78]{dfm1}.
However, the notion of Green function has not been properly introduced here, and we do not want to do so, since the only property of the Green function that we really need is the fact that there exists a solution to $L^*u=0$ in $B(x,r) \cap \Omega$ that satisfies the bounds \eqref{tcp18}.
\end{enumerate}
\end{definition}

The combination of \eqref{Holder} and \eqref{Harnack} gives the existence of $c_2>0$ such that, for any $x\in\partial\Omega$, any $r>0$, any Corkscrew point associated to $(x,r)$, and any Borel set $E\supset \Delta(x,r)$, one has
\begin{equation} \label{NonDeg}
\omega_L^X(E) \geq c_2,
\end{equation}
where $c_2>0$ is independent of $x$, $r$, $X$ and $E$,  and depends on $L$ only via $C_L$.  Indeed, \eqref{Holder} gives that 
\[\begin{split}
\omega_L^{X'}(\partial \Omega \sm E) \leq C \Big( \frac{\delta(X')}{\dist(X',E)} \Big)^\alpha \leq C \Big( \frac{|X'-x|}{|X'-x| - r} \Big)^\alpha \leq  \frac12 
\end{split}\]
as long as $|X'-x| \leq c'r$ with a constant $c'$ that depends only on $C$ and $\alpha$. So if $|X'-x| \leq c'r$  but is still a Corkscrew point associated to $(x,c'r)$, since $\omega^{X'}$ is a probability measure, we have $\omega^{X'}(E) \geq \frac12$. We conclude \eqref{NonDeg} by linking $X'$ and $X$ by a (uniformly finite) Harnack chain of balls and using the Harnack inequality (Lemma \ref{Harnack}) on each of the balls in the chain. 
 
\subsection{Examples of PDE friendly domains}\label{sec.examples}  Let us first state precisely some definitions of boundary  conditions which we have alluded to in previous sections.

\begin{definition}[Ahlfors-David regular set]\label{def.adr} Fix $d\in(0,n-1]$. We say that $\Gamma \subset \R^n$ is a \emph{$d$-Ahlfors-David regular set} (or $d-$ADR) if there exists $C_d > 0$ and a measure $\sigma$ on $\Gamma$ such that 
\begin{equation} \label{defAR}
C_d^{-1} r^d \leq \sigma(B(x,r) \cap \Gamma) \leq C_dr^d \quad \text{ for each } x\in \Gamma, \qquad 0 < r \leq \diam \Gamma.
\end{equation}
If \eqref{defAR} is verified, maybe to the price of taking a larger $C_d$, we can always choose $\sigma$ to be the $d$-dimensional Hausdorff measure on $\Gamma$. 
\end{definition}

\begin{definition}[Capacity and capacity density condition]\label{def.cdc} Given an open set $D\subset\bb R^n$, $n\geq2$, and a compact set $K\subset D$, we define the \emph{capacity} of $K$ relative to $D$ as
\[
\operatorname{Cap}_2(K,D)=\inf\Big\{\dint_D|\nabla v(X)|^2\,dX: v\in C_c^{\infty}(D),~ v(x)\geq1 \text{ on }K\Big\}.
\]
An open set $\Omega$ is said to satisfy the \emph{capacity density condition} (CDC) if there exists a uniform constant $c_1>0$ such that
\[
\frac{\operatorname{Cap}_2(\overline{B(x,r)}\backslash\Omega, B(x,2r))}{\operatorname{Cap}_2(\overline{B(x,r)}, B(x,2r))}\geq c_1,
\]
for all $x\in\partial\Omega$ and $0<r<\diam(\partial\Omega)$.
\end{definition}

We now describe several   examples of PDE friendly domains $(\Omega,m)$.
 
\begin{enumerate}[(i)]
\item\label{ex1} $\Omega$ is a 1-sided NTA domain satisfying the capacity density condition, and $dm=dX$. The elliptic theory for these operators may be found in \cite{hmtbook}, and see \cite{ahmt} for an (additive) perturbation theory in this context. These domains include, in particular, the 1-sided chord-arc (that is, 1-sided NTA and $(n-1)-$ADR) domains, and the $(n-1)-$ADR domains with uniformly rectifiable boundaries.
\item\label{ex2} The case where the boundary is low dimensional; fix $d\in(0,n-1)$ and assume that $\Gamma\subset\bb R^n$ is a $d-$ADR closed set. Set $\Omega=\bb R^n\backslash\Gamma$ and $dm=\delta(X)^{d+1-n}\,dX$. In this situation, the Harnack Chain condition and the Corkscrew point condition are always true, and the elliptic theory was constructed in \cite{dfm1}. An additive perturbation theory was written in \cite{mp20}, for $d\geq1$. Our perturbation theory works for $d\in(0,1)$ as well.
\item\label{ex3} Given $(\Omega,m)$ as in the previous case (and assume that $d\geq1$), given a family $\m F$ of pairwise disjoint ``dyadic cubes'' (these are the David-Christ cubes; see Section \ref{S3} for the definition) on the boundary $\Gamma = \partial \Omega$,  we may construct a sawtooth domain  $(\Omega_{\m F},m_{\m F})$ (with $m_{\m F}=m|_{\Omega_{\m F}}$) that ``hides'' $\m F$, via the procedure in \cite{hm2} using Whitney cubes; see Sections 3 and 4 of \cite{mp20} for the details. The boundary $\partial\Omega_{\m F}$ has pieces of dimension $d$, and other pieces of dimension $n-1$, and thus is mixed-dimensional. In Section 5 of \cite{mp20}, it is shown that we may define a Borel regular measure $\sigma_{\star}$ on $\partial\Omega_{\m F}$ (a ``surface measure'') so that the triple $(\Omega_{\m F}, m_{\m F}, \sigma_{\star})$ satisfies the assumptions (H1)-(H6) of the recent axiomatic elliptic theory in \cite{dfm20}. Briefly, the assumptions (H1)-(H6) include the same conditions on $m_{\m F}$ that we have placed on $m$ in the previous section, the interior Corkscrew point condition and the Harnack Chain condition in $\Omega_{\m F}$, a doubling property of $\sigma_{\star}$, and a slow growth condition on $m_{\m F}$ with respect to $\sigma_{\star}$. With these assumptions verified, the elliptic theory of \cite{dfm20} gives us a sufficiently robust elliptic theory on these mixed-dimensional sawtooth domains. Thus $(\Omega_{\m F}, m_{\m F})$ is a PDE friendly domain.
\item\label{ex4} For that matter, any triple $(\Omega, m, \mu)$ (with $\mu$ a measure on $\partial\Omega$)  satisfying (H1)-(H6) of \cite{dfm20} is a PDE friendly domain. This includes the 1-sided chord-arc domains, the domains with low-dimensional boundaries as in \ref{ex2}, and the mixed-dimensional sawtooth domains as in \ref{ex3}, but there are many more examples, including some domains with boundaries having atoms, and   $t$-independent degenerate  operators on $\R^{n}_+$ that can be written as $L=-\div A(x) \nabla$ with a matrix $A(x)$ that lies in the Muckenhoupt class $A_2(\R^{n-1})$; see Section 3 of \cite{dfm20} for more details and examples.  
\end{enumerate}

\subsection{Local theory} \label{SS2.3}

An important remark is that, even if our main results (Theorems \ref{Main1}, \ref{Main2}, and \ref{Main3}) are stated in their global form, all our theory is intrinsically local. The local versions of our theorems are Corollary \ref{Main11}, Lemma \ref{S<N},   and Lemma \ref{lm.jn2}.

A consequence of the locality of our results and proofs is the fact that, if we are only interested in local results, then we only need to assume a local version of the fact that $(\Omega,m)$ is PDE friendly.

Given $x_0\in \partial \Omega$ and $r_0>0$, we say that $(\Omega,m)$ is \emph{locally PDE friendly} in $B_0:=B(x_0,r_0)$ if the following assumptions hold. 
\begin{enumerate}[(a)]
\item\label{item.a} In Definition \ref{def1NTA}, we only require the existence of Corkscrew points associated to $y\in \frac34B_0$ and $r<r_0/2$. We want Harnack chains between $X$ and $Y$ only when $X,Y\in \frac34B_0 \cap \Omega$. But the Harnack chains might go out of $B_0$, which complicate the definitions. In the sequel, we write $T_0$ for the union of $B_0 \cap \Omega$ with all the Harnack chains linking two points $X,Y\in \frac34B_0 \cap \Omega$.
\item The measure $m$ needs to be defined only on $T_0$ (as given in  \ref{item.a}). We only need the doubling property \eqref{mdoubling} when $B(X,2r) \cap \Omega \subset T_0$, and we only need the Poincar\'e inequality \eqref{defPoincare} when $B \subset T_0$ and $2B \subset \Omega$. 
\item  We need a notion of elliptic measure associated to our elliptic operator $L$. We want a collection $\{\omega^X_L\}_{X \in \Omega \cap B_0}$ of measures on $\Delta_0 := B_0 \cap \partial \Omega$ that satisfy $c_2 \leq \omega^X_L(\Delta_0) \leq 1$ for any $X \in T_0$, and for which $u_E(X):= \omega^X_L(E)$ are solutions to $Lu=0$ in $T_0$ whenever $E\subset \Delta_0$ are Borel sets.
\item On the ``elliptic measure'' defined on the previous point, we ask for the doubling property \eqref{dphm1} only when $\Delta(x,2r) \subset \Delta_0$ and $X\in T_0$, for the change of pole property \eqref{CP18} only when the two poles $X,Y$ are in $T_0$, for the H\"older continuity \eqref{Holder} only when $X\in B_0 \cap \Omega$, and for the comparison with a ``Green function'' \eqref{tcp18} only when $X \in B_0 \cap \Omega$ and $\Delta(y,s) \subset \Delta_0$. 
\end{enumerate}
Under the above assumptions, local forms of our results hold, and we have that if $L_1$ and $L_0$ differs by a Carleson perturbation in $T_0$ (the Carleson perturbations are defined with a local analogue of $KCM$ or $KCM_{\sup}$), then the $A_\infty$ absolute continuity of the harmonic measure is preserved from $L_0$ to $L_1$ in $\frac12B_0 \cap \partial \Omega$. Let us now give an example.  

Assume that $\Omega \subset \R^n$ is an unbounded connected component of an  $(n-1)$-Ahlfors regular set (we call $\sigma$ the Ahlfors regular measure), and assume that $\Omega$ is 1-sided NTA. For instance, $\Omega$ can be $\R^n \setminus B(0,1)$. We choose then $m=\mathcal L^n$ to be the Lebesgue measure on $\Omega$. Then Definitions \ref{defm} and \ref{def1NTA} are verified. So in order to be able to apply our theory, we need to check whether $(\Omega,\mathcal L^n)$ is PDE friendly, and in particular, whether it has a nice elliptic theory. 

However, the harmonic measure $\{\omega^X\}$ associated to this domain - which is the elliptic measure associated to the Laplacian - is not a probability measure, because the Brownian motion has a non-zero probability to escape at infinity. Even worse, we have that $\omega^X(\partial \Omega)$ tends to 0 as $|X| \to \infty$. So we deduce that such $(\Omega,\mathcal L^n)$ are not PDE friendly.  On the other hand, by considering local perspectives, our theory can still be applied. 

The first local perspective is to stick to local estimates. Let $x_0$ be any point in $\partial \Omega$ and for instance $r_0= 100 \diam(\partial \Omega)$. Then one can check that $(\Omega,\mathcal L^n)$ is locally PDE friendly  in $B_0 := B(x_0,r_0)$. Take $L_0$ a uniformly elliptic operator and we assume that its elliptic measure satisfies $\omega_0 \in A_\infty(\sigma, \Delta_0)$, that is, for any $\xi>0$, there exists $\zeta>0$ such that for any boundary ball $\Delta:= \Delta(x,r) \subset \Delta_0:= \Delta(x_0,r_0)$ and any Borel set $E \subset \Delta$, we have
\begin{equation} \label{defAinftyloc}
    \frac{\omega_0^X(E)}{\omega_0^X (\Delta)} < \zeta\quad \text{ implies }\quad \frac{\sigma(E)}{\sigma(\Delta)} < \xi \qquad \text{ for } X\in [\Omega \setminus B(x,2r)] \cap B_0. 
\end{equation}
We say that $f\in KCM_{\Delta_0}(\sigma)$ if there exists $M>0$ such that for $x\in \partial \Omega$ \emph{and $r< r_0$}, we have 
\begin{equation}
\int_{\Delta(x,r)} |\n A^r(f)(y)|^2 \, d\sigma(y) \leq M \sigma(\Delta(x,r)).
\end{equation}
Choose $L_1 = - \div [b(\mathcal A_0 + \mathcal C + \mathcal T) \nabla ]$ a Carleson perturbation of $L_0$ in $B_0$, that is
\[|\mathcal C| \in KCM_{\sup,\Delta_0}(\sigma) \quad \text{ and } \quad \delta|\nabla b| + \delta |\div \mathcal T| \in KCM_{\Delta_0}(\sigma)\]
By applying Lemma \ref{S<N} and Corollary \ref{Main11}, we obtain that $\omega_1 \in A_\infty(\sigma,\Delta_0)$.

An alternative perspective  is to change $(\Omega,\mathcal L^n)$ so that it becomes PDE friendly. We do not want to change $\Omega$, and we want to keep $m=\mathcal L^n$ when we are close to $\partial \Omega$. We need to pay a price on the part of the measure $m$ (and thus the elliptic operators under consideration) which is far from the boundary. We usually do not mind to change the operator far from the boundary, because we can use a comparison principle. We take $m = \mathcal L^n$ on $B_0 \cap \Omega$ and $dm = \delta^{1-n} dx$ on $\R^n \setminus B_0$, where $\delta:=\dist(\cdot,\partial\Omega)$. The reader can then check that the triple $(\Omega,m,\sigma)$ satisfies the assumptions (H1)--(H6) in \cite{dfm20} and hence $(\Omega,m)$ is PDE friendly.

\addtocontents{toc}{\protect\setcounter{tocdepth}{1}}

\section{Theory of $A_\infty$-weights.} \label{S3}

In this section, we gather  the properties of $\partial \Omega$. We do not really need to know that $\Omega$ is PDE friendly, because the incoming results hold on $\partial \Omega$ as a set, except of course when we ultimately apply the theory for elliptic measures in the proof of Theorem \ref{Main3}. The measure $m$ will be mentioned through its appearance in the definition of the area integral $\n A^r$ and hence in the definition of Carleson measure, but the reader can check that none of the properties of $m$ matter.

As before, $\Delta(x,r)$ stands for the boundary ball $B(x,r) \cap \partial \Omega$.  The results on this section will be stated in a ``local'' form, so that they can be applied when $\sigma$ is either a single doubling Borel measures or an elliptic measure (i.e. a collection of measures). 

Let $\Delta_0 = \Delta(x_0,r_0)$ be a boundary ball with $r_0\in(0,\diam\dr \Omega)$. We say that $\sigma$ is doubling in $\Delta_0$ if
\begin{equation} \label{locdoubling}\nonumber
\sigma(\Delta(x,2r) \cap \Delta_0) \leq C_\sigma \sigma(\Delta(x,r) \cap \Delta_0) \qquad \text{ for } x\in \partial \Omega, \, r>0,
\end{equation}
and we say that $\sigma$ is locally doubling if $\sigma$ is doubling in all the boundary balls $\Delta$ (but the constant $C_\sigma$ might depend  on $\Delta$). We say that $f\in KCM_{\Delta_0}(\sigma,M)$ if for any boundary ball $\Delta = \Delta(x,r) \subset \Delta_0$, we have that  $\int_{\Delta} |\n A^r(f)(y)|^2 d\sigma(y) \leq M\sigma(\Delta)$. 
At last, we say that $\sigma_0$ is $A_\infty$-absolutely continuous with respect to $\sigma_0$ on $\Delta_0$ - or $\sigma_1 \in A_\infty(\sigma_0,\Delta_0)$ for short - if Definition \ref{def.ainfty} holds when we assume that $\Delta \subset \Delta_0$, that is, for any $\xi>0$, there exists $\zeta>0$ such that for any boundary ball $\Delta \subset \Delta_0$ and any Borel set $E \subset \Delta$, we have
\begin{equation}
    \frac{\sigma_1(E)}{\sigma_1(\Delta)} < \zeta\quad \text{ implies }\quad \frac{\sigma_0(E)}{\sigma_0(\Delta)} < \xi. 
\end{equation}

We begin the section with some preliminary work on the functional $\n A^r$ introduced in \eqref{defA}. For $\alpha \geq 2$ and $x\in\partial\Omega$, define the cone with larger aperture 
\[\gamma_\alpha^r(x) := \{X\in \Omega, \, |X-x| \leq \alpha \delta(X) \leq \alpha r\}\]
and corresponding area integral
\begin{equation} \label{defAalpha}\nonumber
\n A^r_\alpha(f)(x) := \Big(\dint_{\gamma^r_\alpha(x)} |f(X)|^2 \frac{dm(X)}{m(B_X)} \Big)^\frac12,\qquad f\in L^2_{\loc}(\Omega,m).
\end{equation}
Our first result compares $\n A^r_\alpha$ and $\n A^r$, and is a classical consequence of Fubini's theorem.
 
\begin{lemma}[Comparison of area integrals with different apertures]\label{lm.wide} Let $\alpha \geq 2$,  $\Delta_0 := \Delta(x_0,r_0)$ be a boundary ball with $r_0>0$,  and let $\sigma$ be a doubling measure on $(2+\alpha)\Delta_0$. For any $\Delta_r \subset \Delta_0$ and any $f\in L^2_{\loc}(\Omega,m)$, we have that
	\begin{equation} \label{eq.wide11}
	\dashint_{\Delta_r}|\n A^r_{\alpha}(f)|^2\,d\sigma \leq C_{\alpha}  \dashint_{(2+\alpha)\Delta_r}|\n A^r(f)|^2\,d\sigma,
	\end{equation}
	where $C_{\alpha}$ depends only on $\alpha$ and   $C_\sigma$. Thus, if $f\in KCM_{(2+\alpha)\Delta_0}(\sigma,M_f)$, then
	\begin{equation}\label{eq.wide1}
	M_{\alpha,f}:= \sup_{\Delta_r \subset \Delta_0} \frac1{\sigma(\Delta_r)}\int_{\Delta_r}|\n A^r_{\alpha}(f)|^2\,d\sigma \leq C_{\alpha} M_f.
	\end{equation}
	
\end{lemma}  

\noindent\emph{Proof.} The bound \eqref{eq.wide1} is a straightforward consequence of (\ref{eq.wide11}), which is the only inequality that we need to prove. Fix a surface ball $\Delta_r =\Delta(x,r)\subset\partial\Omega$, and write $T_{\Delta_r}$ for $\bigcup_{y\in \Delta_r} \gamma^r_\alpha(y)$. 
Fubini's lemma entails that
\begin{equation} \label{bdAralpha}
\begin{split}
\int_{\Delta_r}|(\n A_{\alpha}^rf)|^2 \,d\sigma &
= \int_{\partial\Omega} \dint_{\Omega}{\1}_{\Delta_r}(y){\1}_{\gamma_{\alpha}^r(y)}(Y)|f(Y)|^2\,\frac{dm(Y)}{m(B_Y)} \, d\sigma(y) \\
& = \dint_{T_{\Delta_r}} |f(Y)|^2 \frac{dm(Y)}{m(B_Y)}  \Big( \int_{\partial \Omega} {\1}_{\Delta_r}(y){\1}_{\gamma_{\alpha}^r(y)}(Y) \, d\sigma(y) \Big).
\end{split}
\end{equation}
However, $Y\in \gamma^r_\alpha(y)$ if and only if $\delta(Y) \leq r$ and $y\in 4\alpha B_Y$. We deduce that 
\[\begin{split}
\int_{\partial \Omega} {\1}_{\Delta_r}(y){\1}_{\gamma_{\alpha}^r(y)}(Y) \, d\sigma(y) & 
= \sigma(\Delta_r\cap4\alpha B_Y)
\end{split}\]
Let $z$ be such that $|Y-z| = \delta(Y) \leq r$. Then the doubling property of $\sigma$ yields that 
\[\sigma(\Delta_r \cap 4\alpha B_Y) \leq \sigma( \Delta(z,2\alpha\delta(Y))) \lesssim  \sigma( \Delta(z,\delta(Y))) \lesssim \sigma(\partial\Omega\cap8B_Y). \]
The last two computations give  that
\[\begin{split}
\int_{\partial \Omega} {\1}_{\Delta_r}(y){\1}_{\gamma_{\alpha}^r(y)}(Y) \, d\sigma(y)  \lesssim \sigma(\partial\Omega\cap8B_Y) = \int_{\partial \Omega} {\1}_{\gamma^r(y)}(Y) \, d\sigma(y).
\end{split}\]
We reinject the last bound in \ref{bdAralpha} to get
\begin{equation} \label{bdAralpha2}
\begin{split}
\int_{\Delta_r}|(\n A_{\alpha}^rf)|^2 \,d\sigma &
\lesssim \dint_{T_{\Delta_r}} |f(Y)|^2 \frac{dm(Y)}{m(B_Y)}  \Big( \int_{\partial \Omega} {\1}_{\gamma^r(y)}(Y) \, d\sigma(y) \Big) \\
& =  \int_{\partial \Omega} \dint_{T_{\Delta_r}} {\1}_{\gamma^r(y)}(Y) |f(Y)|^2 \frac{dm(Y)}{m(B_Y)} d\sigma(y).
\end{split}
\end{equation}
Observe that for each $y\in\partial\Omega$, $T_{\Delta_r} \cap \gamma^r(y) \neq\varnothing$ precisely when we can find $Y \in \Omega$  and $z\in \Delta_r$ such that $Y \in \gamma_\alpha^r(z) \cap \gamma^r(y)$, and hence
\[|y-x|\leq|y-Y|+|Y-z|+|z-x|< \delta(Y) + \alpha \delta(Y) +r \leq (2+\alpha)r.\]
Consequently, \eqref{bdAralpha2} becomes
\begin{multline}\nonumber
\int_{\Delta_r}|(\n A_{\alpha}^rf)|^2 \,d\sigma  
\lesssim \dint_{T_{\Delta_r}} |f(Y)|^2 \frac{dm(Y)}{m(B_Y)}  \Big( \int_{\partial \Omega} {\1}_{\gamma^r(y)}(Y) \, d\sigma(y) \Big) \\
=  \int_{(2+\alpha)\Delta_r} \dint_{\gamma^r(y)} |f(Y)|^2 \frac{dm(Y)}{m(B_Y)} d\sigma(y) =  \int_{(2+\alpha)\Delta_r} |\n A^r(f)|^2 d\sigma.
\end{multline}
We conclude that
\[ \frac1{\sigma(\Delta_r)}\int_{\Delta_r}|\n A_{\alpha}^r(f)|^2 \,d\sigma \leq C_\alpha \frac1{\sigma((2+\alpha)\Delta_r)} \int_{(2+\alpha)\Delta_r} |\n A^r(f)|^2 d\sigma\]
by using the doubling property of $\sigma$ again. The lemma follows.
\ep

We use the dyadic decomposition of $\partial\Omega$ by Christ, which is a consequence of the metric structure of $\partial \Omega$ induced by $\R^n$. 

\begin{lemma}[Dyadic cubes for a space of homogeneous type \cite{christ}]\label{lm.dyadiccubes} There exists a universal constant $a_0$ such that for each $k\in\bb Z$, there is a collection of sets \emph{(the sets are called ``dyadic cubes'')}
	\[\bb D^k=\bb D^k(\partial\Omega):=\{Q_j^k\subset\partial\Omega\,:\,j\in\n J^k\},\]
	satisfying the following properties.
	\begin{enumerate}[(i)]
		\item\label{item.grid} $\partial \Omega =  \bigcup_{j\in\n J^k}Q_j^k$ for each $k\in\bb Z$, 
		\item\label{item.nonoverlap} If $m\geq k$ then either $Q_i^m\subset Q_j^k$ or $Q_i^m\cap Q_j^k=\varnothing$.
		\item\label{item.generations} For each pair $(j,k)$ and each $m<k$, there is a unique $i\in\n J^m$ such that $Q_j^k\subset Q_i^m$. When $m=k-1$, we call $Q_i^m$ the \emph{dyadic parent} of $Q_j^k$, and we say that $Q_j^k$ is a \emph{dyadic child} of $Q_i^m$.
		\item\label{item.diam} diam $Q_j^k < 2^{-k}$.
		\item\label{item.inscribe} Each $Q_j^k$ contains some surface ball $\Delta(x_j^k,a_02^{-k})=B(x_j^k,a_02^{-k})\cap\partial\Omega$.
	\end{enumerate}
\end{lemma}

\begin{remark}
	The result of Christ assumes a doubling measure on $\partial \Omega$. 
	However, we note that the collections $\bb D^k$ themselves are defined only through the quasi-metric structure of the space, and with no dependence on a doubling measure.
	
	The underlying doubling measure in the result of Christ is used to prove an extra property that imposes a thin boundary (in a quantitative way) on dyadic cubes. Since we do not need thin boundaries for our proofs, and since the result is a bit technical, we avoided to write the full statement. But note that Christ proved, in particular, that for any locally doubling measure $\sigma$ and any dyadic cube $Q_j^k$, we have that $\sigma (\partial Q_j^k) = 0$.
	
	At last, Christ's result only provides the existence of a small $\tau>0$ such that the properties \ref{item.diam} and \ref{item.inscribe} holds for $\tau$ instead of the $1/2$ of our statement. However, it is easy to see that by repeating the collection $\D^k$ over several generations and by taking a smaller $a_0$, it is always possible to take $\tau = 1/2$.
\end{remark}

We shall denote by $\bb D=\bb D(\partial\Omega)$ the collection of all relevant $Q_j^k$; that is,
\[\bb D=\bb D(\partial\Omega):=\bigcup_{k\in\bb Z}\bb D^k(\partial\Omega).\]
Henceforth, we refer to the elements of $\bb D$ as \emph{dyadic cubes}, or \emph{cubes}. For $Q\in\bb D$, we write 
\[\bb D_Q:=\{Q'\in\bb D\,:\,Q'\subseteq Q\} \ \text{  and } \ \bb D_Q^k=\bb D^k(\partial\Omega)\cap\bb D_Q.\]

{\bf Note carefully} that if $Q_i^{k+1}$ is the dyadic parent of $Q_j^k$, then it is possible that, \emph{as sets}, $Q_i^{k+1}=Q_j^k$. In fact, some dyadic cubes may consist of single points (\emph{atoms}), that is a dyadic cube can be equal (as sets) to all of its dyadic descendants. Even if there are no atoms, a dyadic cube could still equal (as sets) an arbitrarily large number of its descendants. Dyadic cubes which are of different generations but are equal as sets, will always be considered distinct. Hence, for $Q\in \mathbb D$, we write $\ell(Q)=2^{-k}$ (the \emph{length} of $Q$) for the {\em only} $k\in \bb Z$ such that $Q\in \mathbb D_k$.

Properties \ref{item.diam} and \ref{item.inscribe} imply that for each  $Q\in\bb D$, there is a point $x_Q\in\partial\Omega$ such that  
\begin{equation}\label{eq.centerofQ} 
\Delta(x_Q,a_0\ell(Q))\subset Q\subseteq \Delta(x_Q,\ell(Q)).
\end{equation}
We call $x_Q$ the \emph{center} of $Q$. 

We redefine our notions using the dyadic cubes instead of the surface balls.  

\begin{definition}[Dyadically doubling measures]\label{def.dyadicdoubling} We say that a Borel measure $\sigma$ on $Q_0\in\bb D$ is \emph{dyadically doubling} in $Q_0$ if $0<\sigma(Q)<\infty$ for every $Q\in\bb D_{Q_0}$ and there exists a constant $C\geq1$ such that $\sigma(Q)\leq C\sigma(Q')$ for every $Q\in\bb D_{Q_0}$ and for every dyadic child $Q'$ of $Q$.
\end{definition}

We let the reader check that if $\sigma$ is a doubling measure in $\Delta_0$ and $Q_0 \subset \Delta_0$, then $\sigma$ is dyadically doubling in $Q_0$.

\begin{definition}[Dyadic $A_{\infty}$ for families of measures]\label{def.ainftyd} Fix $Q_0\in\bb D$. If $\sigma_0$ and $\sigma_1$ are two doubling  measures on $Q_0$, then we say that  $\sigma_1 \in A_\infty^{\dyadic}(\sigma_0,Q)$  if, for any $\xi>0$, there exists $\zeta>0$ such that any $Q\in\bb D_{Q_0}$,  and any Borel set $E \subset Q$, we have that
\begin{equation} \label{defAinftyQ}\nonumber
\frac{\sigma_1(E)}{\sigma_1(Q)} < \zeta \text{ implies } \frac{\sigma_0(E)}{\sigma_0(Q)} < \xi. 
\end{equation}
\end{definition}

We define the truncated area integrals adapted to a dyadic cube $Q \in \bb D$ as
\begin{equation} \label{defAQ}\nonumber
\n A^Q := \n A^{\ell(Q)} \ \text{ and } \n A^Q_\alpha := \n A^{\ell(Q)}_\alpha. 
\end{equation}

\begin{definition}[Dyadic Carleson measure condition] Fix $Q_0\in\bb D$. If $\sigma$ is a doubling measure on $Q_0$, we say that a function $f\in L^2_{\loc}(\Omega,m)$ satisfies the \emph{dyadic $\sigma$-Carleson measure condition on $Q_0$}, written $f\in KCM_{Q_0}(\sigma)$, if there exists $M>0$ such that 
\begin{equation} \label{defCMQ}\nonumber
\int_{Q} |(\n A^{Q}f)(y)|^2 \, d\sigma(y) \leq M \sigma(Q),\qquad\text{for each }Q\in\bb D_{Q_0}.
\end{equation}
We write   $f\in KCM_{Q_0}(\sigma,M)$ if we want to refer to the constant in the above inequality.
\end{definition}

Due to \eqref{eq.centerofQ}, one can see that the dyadic versions of the doubling measure property, the $A_\infty$ absolute continuity, and the Carleson measure condition are {\em a priori} a bit weaker than the general version on balls. However, we can recover the general statement on balls from the dyadic statement, and this is essentially because of the next lemma, which is a slightly refined variant of Lemma 19 in \cite{christ}.

\begin{lemma}[Covering lemma for boundary balls  \cite{christ}]\label{lm.cover} Fix a boundary ball $\Delta:=\Delta(x,r)$, an integer $k\in\bb Z$ such that $a_02^{-k}>r$, and let $\sigma$ be a doubling measure in $\Delta(x,2^{4-k})$.  Then there exists $N\in\bb N$ (depending only on $C_\sigma$ and not on $x,r,k$) such that   there exist  at most $N$ cubes $Q_1^k, \dots , Q_{N_\Delta}^k$ of $\mathbb D^k$ that intersect $\Delta$. 

Consequently, the property \ref{item.grid} of the dyadic decomposition entails that $\Delta \subset  \bigcup_{j=1 }^{N_\Delta} Q_j.$
\end{lemma}

\bp  Let $\Delta:= \Delta (x,r)$ and $k$ be as in the lemma, and let $\{Q_j^k\}_{j\in J}$ be the collection of dyadic cubes in $\bb D^k$ that intersect $\Delta$. Since the number of dyadic cubes is countable, we can identify $J$ to $\{0,\dots, N_\Delta-1\}$ or $\bb N_0$. Due to \eqref{eq.centerofQ}, for each $j\in J$, the center $x_j$ of $Q_j^k$ necessarily satisfies $|x-x_j| \leq r + 2^{-k} \leq 2^{1-k}$, and hence $|x_j - x_0| \leq 2^{2-k}$. We deduce, again thanks to \eqref{eq.centerofQ}, that 
\[ \Delta(x_j,a_02^{-k}) \subset Q_j^k \subset \Delta(x_0,2^{3-k}) \subset \Delta(x_j,2^{3-k}) \qquad \text{ for } j\in J.\]
The doubling property of $\sigma$ entails that the smallest and the biggest sets in the inclusion above have similar measure, hence we also have that $C_{\sigma}'\sigma(Q_j^k)\geq\sigma(\Delta(x_0,2^{3-k}))$ with $C_{\sigma}'$ depending only on the doubling constant of $\sigma$ on $\Delta(x,2^{4-k})$. We conclude that 
\[ C'_\sigma\sigma(Q_j^k) \geq  \sigma(\Delta(x_0,2^{3-k})) \geq \sigma\Big( \bigcup_{j\in J} Q^k_j \Big) = \sum_{i\in J} \sigma(Q_i^k) \qquad \text{ for } j\in J,\]
which means that the cardinality of $J$ is finite and bounded by $C'_\sigma$, as desired. \ep
  
Let us state a local equivalence of the $A_{\infty}$ conditions studied in this article. 

\begin{proposition}[Local interplay of $A_{\infty}$ and $A_{\infty}^{\dyadic}$]\label{prop.ainfty} Let $\sigma_0$ and $\sigma_1$ be two locally doubling Borel measures on $\partial\Omega$. The following statements hold.
\begin{enumerate}[(a)]
	\item Fix $\Delta = \Delta(x,r) \subset\partial\Omega$ and $k\in \bb Z$ such that $a_02^{-k} > r$. If $\sigma_1\in A_{\infty}^{\dyadic}(\sigma_0,Q_j^k)$ for each $Q_j^k \in \bb D^k$ that intersects $\Delta$, then $\sigma_1\in A_{\infty}(\sigma_0,\Delta)$.
	\item Fix $Q\in\bb D(\partial\Omega)$. If for some $r> a_0\ell(Q)$ there exists a cover of $Q$ by a family $\{\Delta_j\}_j$ of surface balls of radius $r$ for which $\sigma_1\in A_{\infty}(\sigma_0,\Delta_j)$ for each $\Delta_j$, then $\sigma_1\in A_{\infty}^{\dyadic}(\sigma_0,Q)$.
	\item If $\sigma_0$ and $\sigma_1$ are both doubling, $\sigma_1\in A_{\infty}^{\dyadic}(\sigma_0)$ if and only if $\sigma_1\in A_{\infty}(\sigma_0)$.
\end{enumerate}
\end{proposition}

\begin{remark} \label{rmk.ainfty}
In (a), the constants in $\sigma_1 \in A_\infty(\sigma_0,\Delta)$ depend only on the doubling constants of $\sigma_0$ and $\sigma_1$ in $\Delta(x,2^{4-k})$, and the constants in $\sigma_1 \in A_\infty^{\dyadic} (\sigma_0,Q^k_j)$.   Of course, a similar property holds for (b). 
\end{remark}

\noindent\emph{Proof.} We prove (a); for the other statements we mention only that the proof of (b) is entirely analogous to that of (a), and (c) follows from (a),  (b), and Remark \ref{rmk.ainfty}. 

Fix $\Delta := \Delta(x,r) \subset\partial\Omega$ and $k\in \bb Z$ such that 
$a_0 2^{-k} >r$. Let $\{Q_j\}_j\subset\bb D^k$ be the collection of cubes in $\mathbb D^k$ that intersect  $\Delta$. Now let $\Delta'=\Delta(x',r')\subset\Delta$ be a surface ball, fix $\xi>0$, and let $k'\in\bb Z$ satisfy $r' < a_0 2^{-k'} \leq 2r'$.
We take $\{Q_j'\}_{j\in J} \subset\bb D^{k'}$ to be the cover for $\Delta'$ afforded by Lemma \ref{lm.cover}, and since $k' \geq k$, it is easy to see that $\sigma_1\in A_{\infty}^{\dyadic}(\sigma_0,Q_j')$.

Let $\zeta>0$ be small to be chosen later, and suppose that $E\subset\Delta'$ is a Borel set that satisfies $\sigma_1(E)< \zeta\sigma_1(\Delta')$, and we want to prove that $\sigma_0(E) < C \xi \sigma_0(\Delta')$ for a constant $C>0$ independent of $\Delta'$ and $E$. For each $j\in J$, we have $\Delta'\cap Q_j' \neq\varnothing$, and therefore
\[\Delta'\subset\Delta(x_{Q_j'},2^{-k'}+2r')\subset\Delta(x_{Q_j'}, 4\ell(Q'_j)).\]
Since $\sigma_1$ is locally doubling, then $\sigma_1(\Delta')\lesssim\sigma_1(\Delta(x_{Q_j'},a_0\ell(Q_j')))\leq\sigma_1(Q_j')$,  and thus 
\[ \sigma_1(E\cap Q_j')\leq C\zeta\sigma_1(Q_j') \qquad\text{for each } j \in J, \]
where $C>0$ depends only of the doubling constant of $\sigma_1$ in $\Delta(x,2^{4-k})$. By the $A_{\infty}$ property on $Q_j'$, there exists $\zeta_j$ small enough  (and independent of $E$ and $\Delta'$) such that $\sigma_0(E\cap Q_j') < \xi \sigma_0(Q_j')$ whenever $\zeta \leq \zeta_j$. We take $\zeta = \min_j\{\zeta_j\}$, which is positive since the number of $Q'_j$ is uniformly bounded, and we obtain
\[ \begin{split}
\sigma_0(E)& = \medsum_j \sigma_0(E\cap Q'_j) \leq \xi \sigma_0(\cup_jQ_j')\leq \xi \sigma_0(\Delta(x',r'+\ell(Q_j'))) \leq C\xi \sigma_0(\Delta'),
\end{split} \]
where we used the doubling property of $\sigma_0$ in $\Delta(x,2^{4-k})$ and  $\ell(Q_j')\lesssim r'$.  \ep

\begin{lemma}[Dyadic cubes as a base for the Carleson measure test]\label{lm.dyadictest} Let $\alpha \geq 2$ and $Q_0\in \bb D$, and let $\Delta_0$ be a boundary ball that contains $\Delta(x_Q,\ell(Q))$ for every $Q\in \bb D_{Q_0}$. Take a doubling   measure $\sigma$ on $(2+\alpha) \Delta_0$, and suppose that $f\in KCM_{(2+\alpha)\Delta_0}(\sigma,M_f)$. Then
\begin{equation}\label{eq.dyadictest}\nonumber
M^{\dyadic}_{\alpha,f} := \sup_{Q\in\bb D_{Q_0}} \frac1{\sigma(Q)} \int_Q|\n A_\alpha^Q(f)|^2\,d\sigma \leq  C M_f,
\end{equation}
where $C>0$ depends only on the doubling constant of $\sigma$.
\end{lemma}

\bp We use \eqref{eq.centerofQ} to change the integration on cubes to integration on balls, and then we conclude using  Lemma \ref{lm.wide}.
\ep

We focus now our efforts on the proof of Theorem \ref{Main3}. We first need a Calder\'on-Zygmund decomposition. Its proof is standard, and is left to the reader.

\begin{lemma}[Calder\'on-Zygmund decomposition]\label{lm.czdecomp} Take $Q_0\in\bb D$ and $\sigma$ a dyadically doubling measure on $Q_0$ with doubling constant $C_{\sigma}$. For any function $f\in L^1(Q_0,\sigma)$ and any level $\lambda > \frac1{\sigma(Q_0)}\int_{Q_0}|f|\,d\sigma$, there exists a collection of maximal and therefore disjoint dyadic cubes $\{Q_j\}_j\subset\bb D_{Q_0}$ such that
	\begin{gather*}
	 f(x)\leq\lambda,\quad\text{for }\sigma-\text{a.e. }x\in  Q_0\setminus \bigcup_jQ_j,\label{eq.cz2}\\ \lambda<\frac1{\sigma(Q_j)}\int_{Q_j}f\,d\sigma\leq C_{\sigma}\lambda \label{eq.cz3}.
	\end{gather*}
\end{lemma}

Our next objective is a John-Nirenberg inequality for Carleson measures. 

\begin{lemma}[John-Nirenberg Lemma for Carleson measures]\label{lm.jn} Let $\Delta_0 \subset \partial \Omega$ be a boundary ball, and let  $\sigma$ be a doubling measure on $30 \Delta_0$ with doubling constant $C_\sigma$. Suppose that $f\in KCM_{30\Delta_0}(\sigma,M_f)$. Then for each boundary ball $\Delta=\Delta(x,r) \subset \Delta_0$, we have that
\begin{equation}\label{eq.jn}
	\sigma\big(\big\{y\in \Delta:|(\n A^r f)(y)|^2>t\big\}\big)\leq Ce^{-ct/M_f}\sigma(\Delta) ,\qquad\text{for }t>0,
\end{equation} 
where $c,C>0$ depend only on $C_\sigma$. 

As a consequence, for any $p\in (0,+\infty)$, we have that
	\begin{equation}\label{eq.jn9}
\Big(\frac1{\sigma(\Delta)}\int_{\Delta}|\n A^r(f)|^p\,d\sigma\Big)^{\frac1p}\leq C_p (M_f)^\frac12,
	\end{equation}
where $C_p$ depends only on $C_\sigma$ and $p$.
\end{lemma}

\noindent\emph{Proof.} The second estimate \eqref{eq.jn9} is a easy consequence of H\"older's inequality (when $p<2$) or \eqref{eq.jn} (when $p>2$). So we only need to prove \eqref{eq.jn}.

We take $f\in KCM_{30\Delta_0}(\sigma,M_f)$ and $\alpha := 4$. Fix $\Delta = \Delta(x,r) \subset \Delta_0$. Let $k \in \bb N$ such that $r < 2^{-k} \leq 2r$, and $\{R_j\}_{j\in J}$ be the collection of dyadic cubes in $\bb D^k$ that intersects $\Delta$. Observe that for $j\in J$, the center $x_j$ of $R_j$ verifies $|x_j - x| \leq 2^{-k} + r \leq 3r$, that is 
\begin{equation} \label{QjinD0}
R_j \subset B(x_j,2^{-k}) \subset 5\Delta \subset 5\Delta_0.
\end{equation}
We can easily check that the above inclusions are also true for every descendant of the $R_j$'s. So for any $R \in \bigcup_j \bb D_{R_j}$, we have $R \subset B(x_{R},\ell(R)) \subset 5\Delta_0$. Lemma \ref{lm.dyadictest} entails that
\begin{equation} \label{eq.jn2}
M^{\dyadic}_{\alpha,f} := \sup_{j\in J} \sup_{R \in \bb D_{R_j}} \frac{1}{\sigma(R)} \int_{R} |\n A^{R}_\alpha(f)|^2 \, d\sigma \leq C' M_f< +\infty,
\end{equation} 
for a $C$ depends only on $C_\sigma$ (recall that $\alpha = 4$, so we have no dependence on $\alpha$).

Fix now $t>0$. By property \ref{item.grid}, the $\{R_j\}_j$ covers $\Delta$, and by \eqref{QjinD0}, the $Q_j$'s stay within $5\Delta$. Those two facts, combined with the fact that $\sigma$ is doubling, entail that the desired estimate \eqref{eq.jn} is a consequence of
\begin{equation}\label{eq.jn1}\nonumber
	\sigma\big(\big\{y\in R_j:|(\n A^{R_j} f)(y)|^2>t\big\}\big)\leq Ce^{-ct/M_f}\sigma(R_j) ,\qquad\text{ for $j\in J$,}
\end{equation} 
where $c,C>0$ depends only on $C_\sigma$.

The index $j$ does not matter anymore, so we drop it and we write $Q_0$ for any of the $R_j$. We also write $M'_f$ for $M_{\alpha,f}^{\dyadic}$ to lighten the notation. The problem is now purely dyadic. Since $\sigma$ is doubling, $\sigma$ is also dyadically doubling with a constant $C'_\sigma$ that depends only on $C_\sigma$. By \eqref{eq.jn2}, we have that
\begin{equation} \label{eq.jn3}
\sup_{Q \in \bb D_{Q_0}} \frac{1}{\sigma(Q)} \int_{Q} |\n A^{Q}_\alpha|^2 \, d\sigma \leq M'_f,
\end{equation} 
and we want to prove that 
\begin{equation}\label{eq.jn4}
	\sigma\big(\big\{y\in Q_0:|(\n A^{Q_0} f)(y)|^2>t\big\}\big)\leq Ce^{-ct/M'_f}\sigma(Q_0).
\end{equation} 
Note that the area integral has different aperture in \eqref{eq.jn3} (big aperture) and \eqref{eq.jn4} (small aperture), and it will become important later in the proof.

Perform the Calder\'on-Zygmund decomposition of the area integral with large aperture $|\n A^{Q_0}_{\alpha}(f)|^2$ on $Q_0$, at height $2M'_{f}$. Since $2M'_{f} >\dashint_{Q_0}|\n A^{Q_0}_{\alpha}(f)|^2\,d\sigma$, according to Lemma \ref{lm.czdecomp}   we may furnish a maximal family $\{Q_{1,j}\}\subset\bb D_{Q_0}$ for which we have
\begin{gather}
\label{eq.czgood}  |(\n A^{Q_0}_{\alpha}(f))(y)|^2\leq2 M_f'~\text{ for }\sigma-a.e.~ y\in Q_0\backslash\cup_jQ_{1,j}, \\ \label{eq.cz4}\nonumber
2M'_f <\frac1{\sigma(Q_{1,j})}\int_{Q_{1,j}}|\n A^{Q_0}_{\alpha}(f)|^2\,d\sigma \leq2 C'_\sigma M'_f. 
\end{gather}
Note that the last line above gives that
\begin{equation} \label{eq.cz5}
\sigma (\cup_jQ_{1,j}) < \frac1{2M'_f} \sum_{j} \int_{Q_{1,j}}|\n A^{Q_0}_{\alpha}(f)|^2\,d\sigma \leq \frac1{2M'_f} \int_{Q_0} |\n A^{Q_0}_{\alpha}(f)|^2\,d\sigma \leq \frac12.
\end{equation}

Let us study the difference of the area integral with small aperture on the cube $Q_{1,j}$.
\[ |\n A^{Q_0}(f)(y)|^2-|\n A^{Q_{1,j}}(f)(y)|^2=\dint_{\gamma^{\ell(Q_0)}(y)\backslash\gamma^{\ell(Q_{1,j})}(y)}|f(Y)|^2\,\frac{dm(Y)}{m(B_Y)},\qquad y\in Q_{1,j}.\]

First, say that $Q_{1,j}'\in\bb D_Q$ is the dyadic parent of $Q_{1,j}$, and let us show that $\sigma(Q_{1,j}'\backslash\cup_kQ_{1,k})\neq0$. Indeed, otherwise there is a (possibly finite) subsequence $Q_{1,k_m}$ such that $Q_{1,j}'=\cup_mQ_{1,k_m}\cup Z$, where $\sigma(Z)=0$. By the dyadic nature
\begin{multline*}
\frac1{\sigma(Q_{1,j})}\int_{Q_j'}|\n A^{Q}_{\alpha}(f)|^2\,d\sigma  =\frac1{\sigma(Q_{1,j}')}\sum_m\int_{Q_{1,k_m}}|\n A^{Q}_{\alpha}(f)|^2\,d\sigma \\ >\frac1{\sigma(Q_{1,j}')}\sum_m 2M_f'\sigma(Q_{1,k_m})=2M'_f,
\end{multline*}
but this is a contradiction to the maximality of the collection $\{Q_{1,j}\}$. The claim is established. Now let $y'\in Q_j'\backslash\cup_kQ_{1,k}$ be arbitrary, and observe that for all $y\in Q_j$, 
\[
\gamma^{Q_0}(y)\backslash\gamma^{Q_{1,j}}(y)\subset\gamma_\alpha^{Q_0}(y'),
\]
where $\gamma_\alpha^{Q_0}$ is the wider cone and $\alpha = 4$. Indeed, if $Y\in\Omega$ belongs to the left-hand side above, then $\ell(Q_{1,j})< \delta(Y) \leq \ell(Q_0)$ for free, and furthermore,
\[
|Y-y'|\leq|Y-y|+|y-y'|\leq2\delta(Y) +\ell(Q'_{1,j})<4\delta(Y).
\]
We have thus deduced that
\[
|(\n A^{Q_0} f)(y)|^2-|(\n A^{Q_{1,j}}f)(y)|^2\leq|(\n A^{Q_0}_\alpha f)(y')|^2\quad\text{for each }y\in Q_{1,j}~\text{ and any }y'\in Q_{1,j}'.
\]
Since we may fix $y'\in Q_{1,j}'\backslash\cup_k Q_{1,k}$ such that (\ref{eq.czgood}) holds at $y'$, we have that
\begin{equation}\label{eq.diff}
|(\n A^{Q_0}f)(y)|^2-|(\n A^{Q_{1,j}}f)(y)|^2\leq 2M'_f,\qquad\text{for each }y\in Q_{1,j}.
\end{equation}

We repeat this process. For each $Q_{1,j}$, we apply the Calder\'on-Zygmund decomposition of $|\n A_{\alpha}^{Q_{1,j}}(f)|^2$ on $Q_{1,j}$, at height $2M_{\alpha}$. Thus there exists a sequence of maximal cubes $\{Q_{2,j}\}$ in $\cup_jQ_{1,j}$ such that
\[
\sigma(\cup_jQ_{2,j})\leq\frac1{2M'_f}\sum_j\int_{Q_{1,j}}|\n A_{\alpha}^{Q_{1,j}}(f)|^2\,d\sigma\leq\frac1{2M'_f}\sum_jM'_f\sigma(Q_{1,j})< 2^{-2} \sigma(Q).
\]
Moreover, on $Q_0 \backslash\cup_jQ_{1,j}$, we have that $|(\n A^{Q_0}f)(y)|^2\leq|(\n A^{Q_0}_{\alpha}f)(y)|^2\leq2M'_f$ for $\sigma-$a.e. $y$; while on $\cup_jQ_{1,j}\backslash\cup_iQ_{2,i}$, thanks to (\ref{eq.diff}), for $\sigma-$a.e. $y$   we have that
\[
|(\n A^{Q_0}f)(y)|^2\leq|(\n A^{Q_0}f)(y)|^2-|(\n A^{Q_{1,j}}f)(y)|^2+|(\n A^{Q_{1,j}}_\alpha f)(y)|^2\leq2M'_f+2M'_f=2(2M'_f).
\]
Consequently, $|(\n A^{Q_0}f)(y)|^2\leq 2(2M'_f)$ $\sigma-$a.e. on $Q_0\backslash\cup_k Q_{2,k}$.

We may now iterate this process. As such, for each integer $k\in\bb N$, there exists a sequence of maximal cubes $\{Q_{k,j}\}$ such that $\sigma(\cup_kQ_{k,j})\leq2^{-k}\sigma(Q)$, and (via an easy telescoping argument)
\[
|(\n A^{Q_0}f)(y)|^2\leq 2kM'_f,\qquad\text{for }\sigma-\text{a.e. }y\in Q\backslash\cup_kQ_{k,j}.
\]
Therefore, we have shown that $\sigma(\{y\in Q_0:|(\n A^{Q_0}f)(y)|^2>2kM'_f\})\leq2^{-k}\sigma(Q_0)$ for each integer $k\geq0$, whence (\ref{eq.jn4}) easily follows.\hfill{$\square$}

We recall here a classical characterization of  $A_\infty$ via reverse H\"older estimates. 

\begin{proposition}[$RH$ characterization of $A_{\infty}$ \cite{gr}] \label{pr.Ainfty}
Let $\sigma_0$ and $\sigma_1$ be two doubling measures on a boundary ball $\Delta \subset \partial \Omega$. Then the following are equivalent:
\begin{enumerate}[(i)]
\item $\sigma_1 \in A_\infty(\sigma_0,\Delta)$,
\item $\sigma_0 \in A_\infty(\sigma_1,\Delta)$,
\item $\sigma_1\ll\sigma_0$  and the Radon-Nikodym derivative  $k:=d\sigma_1/d\sigma_0$ satisfies a reverse H\"older bound on $(\Delta,\sigma_0)$. More precisely, there exists $q>1$ and $C>0$ such that
\begin{equation} \label{ReverseHolder}
\Big( \frac1{\sigma_0(\Delta')} \int_{\Delta'} k^q \, d\sigma_0 \Big)^\frac1q \leq C \frac1{\sigma_0(\Delta')} \int_{\Delta'} k \, d\sigma_0 \ \text{ for any boundary ball } \Delta' \subset \Delta.
\end{equation}
If $k$ satisfies (\ref{ReverseHolder}), we say that $k\in RH_q(\Delta,\sigma_0)$.
\end{enumerate}
\end{proposition}

The only time when we need the powerful characterization of $A_\infty$  given above is to prove the following transitivity of Carleson measures.

\begin{lemma}[Local $A_\infty$ implies the transference of the Carleson measure condition] \label{lm.jn2}
Let $\Delta \subset \partial \Omega$ be a boundary ball, and let $\sigma_0,\sigma_1$ be two doubling measures on $30 \Delta$. If $\sigma_1 \in A_\infty(\sigma_0,\Delta)$, then for each $f\in L^2_{\loc}(\Omega,m)$, 
\[
\text{if}\qquad f \in KCM_{30\Delta}(\sigma_0,M_f),\qquad\text{then}\qquad f \in KCM_{\Delta}(\sigma_1,CM_f),
\]
where $C>0$ depends only on the   doubling constant of $\sigma_0$ and the constants $C,q$ in the characterization of $\sigma_1 \in A_\infty(\sigma_0,\Delta)$ given in Proposition \ref{pr.Ainfty}.
\end{lemma}

\bp Let $\Delta$, $\sigma_0$, and $\sigma_1$ be as in the assumption of the lemma, fix $f \in KCM_{30\Delta}(\sigma_0,M_f)$ and $\Delta' \subset \Delta$. We want to prove that $\frac1{\sigma_1(\Delta')}\int_{\Delta'}|\n A^r(f)|^2\,d\sigma_1 \leq C M_f$.  Since $\sigma_1 \in A_\infty(\sigma_0,\Delta')$, writing   $k = d\sigma_1/d\sigma_0$ and H\"older inequality gives that 
\begin{multline*}
\dashint_{\Delta'}|\n A^r(f)|^2\,d\sigma_1   = \frac{\sigma_0(\Delta')}{\sigma_1(\Delta')} \dashint_{\Delta'}|\n A^r(f)|^2\,k \, d\sigma_0 \\
  \leq \frac{\sigma_0(\Delta')}{\sigma_1(\Delta')} \Big(\dashint_{\Delta'} k^q \,d\sigma_0 \Big)^\frac1q \Big(\,\dashint_{\Delta'}  |\n A^r(f)|^{2p} \,d\sigma_0 \Big)^{\frac1p}
\end{multline*}
where $q>1$ is the parameter given by Proposition \ref{pr.Ainfty} and $\frac1p+\frac1q = 1$. Using \eqref{ReverseHolder} and \eqref{eq.jn9} allows us to deduce
\[\begin{split}
\dashint_{\Delta'}|\n A^r(f)|^2\,d\sigma_1 & \lesssim  \frac{\sigma_0(\Delta')}{\sigma_1(\Delta')} \Big(\dashint_{\Delta'} k \,d\sigma_0\Big) M_f = M_f. 
\end{split}\]
The lemma follows. \ep

\noindent {\em Proof of Theorem \ref{Main3}.}   We shall only consider the case where $\mu_0 = \sigma_0$ is a doubling measure and $\mu_1 = \{\omega_1^X\}_{X\in \Omega}$ is an elliptic measure, and we shall only prove the implication 
\[f \in KCM(\sigma_0) \implies f \in KCM(\omega_1),\qquad\text{for each }f\in L^2_{\loc}(\Omega,m).\]
All the other situations are analogous to this one with obvious modifications. 

So take $f\in L^2_{\loc}(\Omega,m)$ that verifies $f\in KCM(\sigma_0,M_f)$. We will show that $f\in KCM(\omega_1,CM_f)$. Thus fix $x\in\dr \Omega$, $r\in(0,\diam\Omega)$, and let $Y\in\Omega$ be a Corkscrew point for $\Delta:=\Delta(x,r)$. There exists $c_1>0$ such that $\delta(Y)\geq 60c_1r$, so $\omega_1^Y$ is doubling on $c_1\Delta=\Delta(x,c_1r)$ by \eqref{dphm1}. However, $\omega_1^Y$ is also doubling on $30\Delta$. Indeed, we can cover $30\Delta$ by a uniformly finite number of small balls $\{\Delta_i = \Delta(x_i,r')\}$ of radius $r'=c_1r/2$ by the Vitali covering lemma, then we pick corkscrew points $Y_i$ associated to $(x_i,r)$, and the same argument yields that $\omega^{Y_i}$ is doubling on $2\Delta_i$. The Harnack chain condition allows us to connect $Y_i$ and $Y$ by Harnack chains, and the Harnack inequality (Lemma \ref{Harnack}) yields that $\omega_1^Y$ is doubling on each ball $2\Delta_i$ and then on $\Delta$.

Of course, by assumption, we also have $f\in KCM_{30\Delta}(\omega_0^Y)$ and that $\omega_1^Y\in A_{\infty}(\omega_0^Y,\Delta)$, so by Lemma \ref{lm.jn2}, we deduce that $f\in KCM_{\Delta}(\omega_1^Y, C'M_f)$, and $C'$ is independent of $\Delta$ and $Y$.

We conclude by the change of pole property \eqref{CP18}, which shows without difficulty that 
\[f \in KCM_\Delta(\omega^Y_1,C'M_f) \implies f \in KCM_\Delta(\omega^X_1, C'' M_f),\qquad\text{ for } X\in \Omega \setminus B(x,2r)\]
for a constant $C''$ independent of $\Delta$ and $X$. The theorem follows.
  \ep

\section{Proof of Theorem \ref{Main1}} \label{S4}

Our proof method is analogous to that of \cite[Theorem 8.9]{dfm2}; see also \cite{kkipt} and \cite{chmt}. In particular, we remark that our method of proof for Theorem \ref{Main1} differs from that of \cite[Theorem 1.1$(a)\implies(b)$]{chmt} in that we do not (and cannot, because it is not true in our more general setting of PDE-friendly domains) use  that every dyadic cube will have a proper descendant after a uniform number of dyadic generations, nor do we use (and cannot use) the largeness of the elliptic measure of the complement of a surface ball. 

We also want to thank Jos\'e-Mar\'ia Martell for pointing out to us that we do not need to assume in our proof that the elliptic  measure is a probability measure, but only that the full measure of the elliptic measure is uniformly bounded from below by a constant $c_2>0$. We changed our proof to match this case.

We will prove the following local result, that implies Theorem \ref{Main1}.

\begin{lemma}[Local $KCM\implies$ local $A_{\infty}$, dyadic version]\label{Main1loc} Let $(\Omega,m)$ be PDE friendly. Let $L= - \div A\nabla$ be an elliptic operator satisfying \eqref{defellipticw} and \eqref{defboundedw}, and construct the elliptic measure $\omega:= \{\omega^X\}_{X\in \Omega}$ as in \eqref{defhm}. 
	
There exists $\alpha \geq 2$ that depends only on the constants in the Corkscrew point condition, the Harnack chain condition, and the H\"older continuity \eqref{Holder} such that the following holds. Fix $Q_0\in\bb D(\partial\Omega)$. If there exists a constant $M>0$ and a dyadically doubling measure $\sigma$ on $Q_0$ such that, for any Borel $E \subset Q_0$, the solution $u_E$ constructed as $u_E(X) := \omega^X(E)$ satisfies
\begin{equation} \label{Main1aloc}
\sup_{Q \in \bb D_{Q_0}} \dashint_Q |\n A_\alpha^Q(\delta \nabla u_E)|^2 d\sigma \leq M,
\end{equation}
then $\omega \in A_\infty^{\dyadic}(\sigma,Q_0)$.
\end{lemma}

The lemma implies 

\begin{corollary}[Local $KCM\implies$ local $A_{\infty}$] \label{Main11}
Let $(\Omega,m)$ be PDE friendly. Let $L$ satisfy \eqref{defellipticw} and \eqref{defboundedw}, and let $\omega:= \{\omega^X\}_{X\in \Omega}$ be the associated elliptic measure. 

There exists $K>0$ that depends only on the same parameters as $\alpha$ in Lemma \ref{Main1loc} such that the following holds. Take $\Delta_0$ to be a boundary ball. If for any Borel $E \subset \Delta_0$, the solution $u_E$ constructed as $u_E(X) := \omega^X(E)$ satisfies $\delta \nabla u_E \in KCM_{K\Delta_0}(\sigma,M)$ for a constant $M>0$ and a doubling measure $\sigma$ on $K\Delta_0$, then $\omega \in A_\infty(\sigma,\Delta_0)$.
\end{corollary}

\noindent {\em Proof of Corollary \ref{Main11} from Lemma \ref{Main1loc}.}
Let $\alpha\geq 2$ as in Lemma \ref{Main1loc} and $K=5(2+\alpha)$. We construct the collection $\{R_j\}_{j\in J}$ of dyadic cubes that covers $\Delta_0$ as in the proof of Lemma \ref{lm.jn}, and the same reasoning as in the proof of Lemma \ref{lm.jn} yields that 
\begin{equation}
M^{\dyadic}_{\alpha} := \sup_{E \subset \Delta_0} \sup_{j\in J} \sup_{R \in \bb D_{R_j}} \frac{1}{\sigma(R)} \int_{R} |\n A^{R}_\alpha(\delta \nabla u_E)|^2 \, d\sigma \leq C' M< +\infty.
\end{equation} 

Lemma \ref{Main1loc} gives then that $\omega \in A_\infty^{\dyadic}(\sigma,R_j)$ for each $j\in J$ and Proposition \ref{prop.ainfty} allows us to recover the non dyadic version $\omega\in A_{\infty}(\sigma,\Delta_0)$. \ep

\noindent {\em Proof of Theorem \ref{Main1}.}
If $\sigma$ is a doubling measure on $\partial \Omega$, then Theorem \ref{Main1} is a straightforward consequence of Corollary \ref{Main11}.

When $\sigma$ is an elliptic measure, Theorem \ref{Main1} is a consequence of Corollary \ref{Main11}, and the properties of the elliptic measure $\sigma$ (doubling property \eqref{dphm1}, change of pole \eqref{CP18}). 
\ep

The rest of the section is devoted to the proof of Lemma \ref{Main1loc}.

\subsection{Step I: Construction of functions with large oscillations on small sets} The first order of business will be to construct the regions over which we will have large oscillations.

\begin{definition}[Good $\eps_0$ cover]\label{def.cover} Fix $Q\in\bb D(\partial\Omega)$ and let $\nu$ be a regular Borel measure on $Q$. Given $\eps_0\in(0,1)$ and a Borel set $E\subset Q$, a \emph{good $\eps_0-$cover} of $E$ with respect to $\nu$, of length $k\in\bb N$, is a collection $\{\m O_{\ell}\}_{\ell=1}^k$ of Borel subsets of $Q$, together with pairwise disjoint families $\m F_{\ell}=\{S_i^{\ell}\}\subset\bb D_{Q}$, such that
\begin{enumerate}[(a)]
	\item \label{item.EinO} $E\subset\m O_k\subset\m O_{k-1}\subset\cdots\subset\m O_2\subset\m O_1\subset\m O_0=Q$,
	\item \label{item.O=S} $\m O_{\ell}=\bigcup_i S_i^{\ell},\qquad 0\leq\ell\leq k$,
	\item \label{item.SO<S} $\nu(\m O_{\ell}\cap S_i^{\ell-1})\leq\eps_0\nu(S_i^{\ell-1})$, for each $S_i^{\ell-1}\in\m F_{\ell-1}$, $1\leq\ell\leq k$.
	\item \label{item.children} for each $S_i^{\ell-1}\in\m F_{\ell-1}$, $1\leq\ell\leq k$, the dyadic cube $S_i^{\ell-1}$ has at least two different children.
\end{enumerate}
\end{definition}

\begin{remark}\label{rm.children} 
The \emph{good $\eps_0-$cover} has already been considered in multiple works, such as \cite{kkipt,dfm2,chmt}. In all those works, the property \ref{item.children} is not stated, but we can actually get this extra assumption for free, as explained in the following lines. First, we can always assume that $S^{\ell}_i$ intersects $E$,  because otherwise we remove each $S^{\ell}_i$ that does not intersect $E$ from the collections $\m F_{\ell}$, and still get the same properties \ref{item.EinO}, \ref{item.O=S}, and \ref{item.SO<S}. With this in hand, $\mathcal O_l \cap S^{\ell-1}_i$ will never be empty, and thus property \ref{item.SO<S} implies that $S^{\ell-1}_i$ cannot be an atom (that is, a set reduced to one point). At last, the cubes $\{S_i^{\ell}\}$ making up the good $\eps_0$-cover are chosen \emph{as sets}, meaning that the generation does not matter, and since $\{S_i^{\ell-1}\}$ are not atoms, we can always choose $S_i^{\ell-1}$ so that its child is not $S_i^{\ell-1}$, meaning that $S_i^{\ell-1}$ possesses at least two children.

\end{remark}

As in \cite{dfm2}, we write  $S_i^\ell$ for the cubes making up $\m O_{\ell}$ so as not to abuse the notation $Q_i^\ell$, which is reserved for a dyadic cube of generation $\ell$. 
Next, we have the fact that we may construct good $\eps_0-$covers. Although the analogous statement in \cite[Lemma 3.5]{chmt} is formally only for the case of $n-$dimensional Ahlfors-David regular sets, a study of their proof reveals no dependence on the Ahlfors regularity \emph{per se}, and their argument extends seamlessly to our setting. See also the remark that follows.

\begin{lemma}[Existence of good $\eps_0-$covers, {\cite[Lemma 3.5]{chmt}}]\label{lm.existgood} Fix $Q\in\bb D(\partial\Omega)$. Let $\nu$ be a doubling measure on $Q$, with dyadic doubling constant $C_{\nu}^{\dyadic}$. For every $0<\eps_0<e^{-1}$, if $E\subset Q$ is a Borel set with $\nu(E)\leq\zeta\nu(Q)$ and $0<\zeta\leq\eps_0^2/(\sqrt2C_{\nu}^{\dyadic})^2$ then $E$ has a good $\eps_0-$cover with respect to $\nu$ of length $k_0=k_0(\zeta,\eps_0)\in\bb N$, $k_0\geq2$, which satisfies 
\[
k_0\gtrsim\frac{\log(\zeta^{-1})}{\log(\eps_0^{-1})}.
\]
In particular, if $\nu(E)=0$, then $E$ has a good $\eps_0-$cover of arbitrary length.
\end{lemma}

\begin{remark}\label{rm.zeroth} The good $\eps_0-$cover constructed in \cite{chmt} does not specify the zeroth cover $\m O_0$; however, it is an easy exercise to check that $\m O_0=Q$ with $\{S_i^0\}=\{Q\}$ can be appended to the cover $\{\m O_{\ell}\}_{\ell=1}^k$ from \cite[Lemma 3.5]{chmt} to produce a good $\eps_0-$cover in our sense of Definition \ref{def.cover}.
\end{remark}

We will eventually show that $\omega\in A_{\infty}^{\dyadic}(\sigma,Q_0)$ (see Definition \ref{def.ainftyd}), but first we need to set the table. Fix $Q\in\bb D_{Q_0}$, let $X_0\in\Omega\backslash B(x_{Q_0},2\ell(Q_0))$. Observe that $\omega^{X_0}$ is a regular Borel measure on $\partial\Omega$ which is dyadically doubling on $Q_0$ by \eqref{dphm1}. Henceforth  we let $0<\eps_0<e^{-1} $ and $0<\zeta<\eps_0^2/(2C_0^2)$ be sufficiently small to be chosen later, and we let $E\subset Q$ be a Borel set such that $\omega^{X_0}(E)\leq\zeta\omega^{X_0}(Q)$. We may apply Lemma \ref{lm.existgood} with $\nu=\omega^{X_0}$ to exhibit a good $\eps_0-$cover for $E$ of length $k\gtrsim\frac{\log(\zeta^{-1})}{\log(\eps_0^{-1})}$ with $k\geq2$. Thus let $\{\m O_{\ell}\}_{\ell=0}^k$ and $\{S_i^{\ell}\}_{\m F_{\ell}}$ be as described in Definition \ref{def.cover}. 

Owing to the property \ref{item.children} of the $\epsilon_0$-cover, for each $S_i^{\ell}$, we let $\widehat S_i^{\ell}$ and $\wt S_i^\ell$ be two different children of $S_i^{\ell}$. Following ideas of \cite{kkipt} and \cite{dpp17}, we set $\widehat{\m O}_{\ell}:=\bigcup_i\widehat S_i^{\ell}\subset\m O_{\ell}$ for each $\ell=0,\ldots,k$. Now, without loss of generality we may take $k$ to be odd, and for each even $\ell$ with $0\leq\ell\leq k-1$, we define
\[
f_{\ell}:={\1}_{\widehat{\m O}_{\ell}},\qquad f_{\ell+1}:=-f_{\ell}{\1}_{\m O_{\ell+1}}=-{\1}_{\widehat{\m O}_{\ell}\cap\m O_{\ell+1}},
\]
so that $f_{\ell}+f_{\ell+1}={\1}_{\widehat{\m O}_{\ell}\backslash\m O_{\ell+1}}$ for $\ell$ even, and write
\begin{equation}\label{eq.f}
f:=\sum_{\ell=0}^kf_{\ell}=\sum_{l=0}^{(k-1)/2}{\1}_{\widehat{\m O}_{2l}\backslash\m O_{2l+1}}={\1}_{\bigcup_{l=0}^{(k-1)/2}\big(\widehat{\m O}_{2l}\backslash\m O_{2l+1}\big)}.
\end{equation} 

\subsection{Step II: The solution with data $f$ exhibits large oscillations on Whitney cubes} Let $u$ solve $Lu=0$ with data $f$ on $\partial\Omega$, and according to (\ref{eq.f}), we have that $u(X)=\omega^X\big(\bigcup_{l=0}^{(k-1)/2}(\widehat{\m O}_{2l}\backslash\m O_{2l+1})\big)$. We shall present two balls, close to one another, over which $u$ oscillates. 

Take any $x\in E$, and $0\leq\ell\leq k$, $\ell$ even. Let $S^{\ell}\in\{S_i^{\ell}\}$ be the unique cube that contains $x$, that possesses (at least) the two children $\wh S^\ell$ and $\wt S^\ell$. We write  $r_\ell$ for $\ell(S^\ell)$, we call $\wh x_\ell$ and $\wt x_\ell$ the centers of $\wh S^\ell$ and $\wt S^\ell$ respectively, and we set $\wh \Delta_\ell := \Delta (\wh x_\ell, a_0r_\ell/2) \subset \wh S^\ell$ and $\wt \Delta_\ell := \Delta (\wt x_\ell, a_0r_\ell/2) \subset \wt S^\ell$.

By the H\"older continuity \eqref{Holder} of the elliptic measure at the boundary, we deduce that there exists $\rho>0$ such that
\begin{equation} \label{zz1}
\omega^{X}(\partial \Omega \sm \wh \Delta_\ell) \leq \frac{c_2}8 \qquad \text{ for } X \in B(\wh x_\ell,\rho r_\ell) \cap \Omega,
\end{equation}
where $c_2$ is the constant from the non-degeneracy  of the elliptic measure (\ref{NonDeg})\footnote{We use the estimate (\ref{NonDeg}) to show that our argument is fundamentally local, and does not depend on the global properties of the elliptic measure; in particular, our argument does not directly use the fact that $\omega(\partial\Omega)=1$.}, and similarly 
\begin{equation} \label{zz2}
\omega^{X}(\partial \Omega \sm \wt \Delta_\ell) \leq \frac{c_2}8 \qquad \text{ for } X \in B(\wt x_\ell,\rho r_\ell) \cap \Omega.
\end{equation}
For the rest of the proof, $\wh X_\ell$ and $\wt X_\ell$ are Corkscrew points associated to respectively $(\wh x_\ell,\rho r_\ell)$ and $(\wt x_\ell,\rho r_\ell)$. That is, for a constant $c$ that depends only on $\rho$, the constant $c_1$ in Definition \ref{def1NTA}, and $c_2$ from (\ref{NonDeg}), we have 
\[B(\wh X_\ell,cr_\ell) \subset B(\wh x_\ell,\rho r_\ell) \cap \Omega\ \text{ and }\ B(\wt X_\ell,cr_\ell) \subset B(\wt x_\ell,\rho r_\ell) \cap \Omega.\] 
So if we set $\wh B_\ell := B(\wh X_\ell,cr_\ell/20)$ and $\wt B_\ell := B(\wt X_\ell,cr_\ell/20)$, the bounds \eqref{zz1}, (\ref{NonDeg}), and \eqref{zz2} entail
\begin{equation} \label{zz3}
\omega^{X}(\wh \Delta_\ell) \geq \frac78c_2 \ \text{ for } X \in \wh B_\ell \qquad \text{ and } \omega^{X}(\partial \Omega \sm \wt \Delta_\ell) \leq \frac{c_2}8 \ \text{ for } X \in \wt B_\ell.
\end{equation}

We want to use the above bounds to estimate $u$ on the balls $\wh B_\ell$ and $\wt B_\ell$. For each $X\in \wh B_\ell$, we have
\begin{equation}\label{eq.ul}
u(X)\geq\omega^{X}\big(\widehat{\m O}_{\ell}\backslash\m O_{\ell+1}\big)\geq\omega^{X}\big(\widehat\Delta_{\ell}\backslash\m O_{\ell+1}\big)=\omega^{X}(\widehat\Delta_{\ell})-\omega^{X}(\widehat\Delta_{\ell}\cap\m O_{\ell+1}).
\end{equation}
and we want to show that the second term of the right-hand side above is small, smaller than $c_2/8$. Observe that
\begin{equation}\label{eq.ul2}
\omega^{X}(\widehat\Delta_{\ell}\cap\m O_{\ell+1})\lesssim \frac{\omega^{X_0}(\widehat\Delta_{\ell}\cap\m O_{\ell+1})}{\omega^{X_0}(\widehat\Delta_{\ell})}\leq\frac{\eps_0\omega^{X_0}(S^{\ell})}{\omega^{X_0}(\widehat\Delta_{\ell})}\lesssim\eps_0\frac{\omega^{X_0}(S^{\ell})}{\omega^{X_0}(S^{\ell})}=\eps_0,
\end{equation}
where we have used the change of pole \eqref{CP18}, then property \ref{item.SO<S} of the good $\epsilon_0$ cover, and at last the doubling property of $\omega^{X_0}$. Therefore, there exists a constant $M$ so that  $\omega^{\widehat X_{\ell}}(\widehat\Delta_{\ell}\cap\m O_{\ell+1})\leq M\eps_0$. If we ask that $\eps_0<c_2/(8M)$, then putting (\ref{eq.ul}), \eqref{zz3}, and (\ref{eq.ul2}) together we may conclude that  
\begin{equation}\label{eq.ul3}
u(X)\geq \frac34c_2,\qquad\text{ for  } X\in\widehat B_\ell.
\end{equation}
Thus we have that $u$ is large on a Whitney region associated to $S^\ell$. Similarly, for $X\in \wt B_\ell$, we have
\begin{equation} \begin{split}\label{eq.decompose}\nonumber
& u(X) =\omega^{X}\big(\medcup_{l=0}^{(k-1)/2}(\widehat{\m O}_{2l}\backslash\m O_{2l+1})\big) \\
& \leq \omega^{X}\big(\partial \Omega \sm \wt \Delta_\ell \big) +\sum_{l=0}^{(k-1)/2} \omega^{X}\big((\widehat{\m O}_{2l}\backslash\m O_{2l+1}) \cap \wt \Delta_\ell \big) \\
& \leq \omega^{X}\big(\partial \Omega \sm \wt \Delta_\ell \big) + \omega^{X}\big( \widehat{\m O}_{\ell} \cap \wt \Delta_\ell \big) + \sum_{2l+1 < \ell} \omega^{X}\big(\wt \Delta_\ell  \backslash\m O_{2l+1} \big) +\sum_{2l > \ell } \omega^{X}\big( \widehat{\m O}_{2l} \cap \wt \Delta_\ell \big).
\end{split}\end{equation}
By construction, $\widehat{\m O}_{\ell} \cap \wt \Delta_\ell = \varnothing$. Notice also that $\wt \Delta_\ell \subset S^\ell \subset \mathcal O_{\ell -1}$, hence $\wt \Delta_\ell \backslash\m O_{2l+1} = \varnothing$ when $2l+1 < \ell$. When $2l > \ell$, using the change of pole \eqref{CP18} and the property \ref{item.SO<S} of the good $\epsilon_0$ cover like in \eqref{eq.ul2}, we obtain for $X\in \wt B_\ell$ that
\begin{equation}\label{eq.ul22}\nonumber
\omega^{X}\big(\widehat{\m O}_{2l} \cap \wt \Delta_\ell\big) \leq \omega^{X}\big(\wt \Delta_\ell \cap \m O_{2l} \big) \lesssim \frac{\omega^{X_0}(\wt \Delta_\ell \cap \m O_{2l})}{\omega^{X_0}(\wt \Delta_\ell)} \leq (\eps_0)^{2l-\ell} \frac{\omega^{X_0}(S^{\ell})}{\omega^{X_0}(\widehat\Delta_{\ell})}\lesssim (\eps_0)^{2l-\ell}.
\end{equation}
Owing to \eqref{zz3} and the observations above, the bound on $u$ when $X\in \wt B_\ell$ becomes  $u(X) \leq \frac{c_2}8 + M' \sum_{2l > \ell} (\eps_0)^{2l-\ell}$ for some $M'$ that is independent of all the important parameters. We choose $\eps_0$ small enough so that $M' \sum_{2l > \ell} (\eps_0)^{2l-\ell} <c_2/8$, and we conclude 
\[u(X)\leq \frac{c_2}4 ,\qquad\text{ for } X \in \wt B_\ell. \]
The last inequality together with (\ref{eq.ul3}) imply the desired large oscillation result. More precisely, if $B\subset\partial\Omega$ is a ball and we write $u_B:=\frac1{m(B)}\iint_Bu\,dm$, then we have that
\begin{equation}\label{eq.osc}
|u_{\widehat B_\ell}-u_{\wt B_\ell}|\geq c_2/2.
\end{equation}

\subsection{Step III: Large oscillations on Whitney regions imply large square function} We now purport to pass from the large oscillation estimate (\ref{eq.osc}) to a pointwise lower bound on the square function.

\subsubsection{A Poincar\'e estimate} We ought to pass from the estimate on the difference over similarly sized balls to an estimate on the gradient, and this can be done via a delicate use of the Poincar\'e inequality. First of all, we recall that the radii of $\widehat B_\ell$ and $\wt B_\ell$ are equivalent to $r_\ell = \ell(S^\ell)$. Moreover, $\widehat B_\ell$, $\wt B_\ell$ are chosen so that both $20 \widehat B_\ell$ and $20 \wt B_\ell$ are subset of $\Omega \cap B(x_{S^\ell},r_\ell)$.
Therefore, we have that $\min\{\delta(\widehat X_\ell),\delta(\wt X_\ell)\}\geq r_\ell /M$ and $|\widehat X_\ell-\wt X_\ell| \leq 2 r_\ell$. The Harnack chain condition from Definition \ref{def1NTA} (and Remark \ref{rm.Harnack}) provides the existence of a  Harnack Chain $\{B_j\}_{j=0}^N=\{B(X_j,\operatorname{rad}(B_j))\}_{j=0}^N$ of balls such that $N$ is a uniformly bounded number (depending only on the allowable constants), $B_0=\widehat B_\ell$, $B_N=\wt B_\ell$, $\delta(X_j) = 20 \operatorname{rad}(B_j)$, and $B_j\cap B_{j+1}\neq\varnothing$ for each $j$ (this last property can be ensured by adding in more balls of the same radius if necessary).
Under this setup, \eqref{eq.osc} becomes
\begin{equation}\label{eq.balls}
\frac{1}2 \leq|u_{B_0}-u_{B_N}|\leq\sum_j|u_{B_j}-u_{B_{j+1}}|\lesssim \sum_j(|u_{B_j}-u_{3B_j}|+|u_{B_{j+1}}-u_{3B_j}|).
\end{equation}
We now assume that $j=j(\ell)$ is the index at which the  maximum in the right-hand side of (\ref{eq.balls}) is taken. Since $B_j\cup B_{j+1}\subset3B_j$, we may estimate
\begin{equation} \begin{split} \label{eq.poincare1}
|u_{B_j}-u_{3B_j}|& =  \fiint_{B_j}|u-u_{3B_j}| \,dm \leq \fiint_{3B_j}|u -u_{3B_j}|\,dm \\ & \lesssim \operatorname{rad}(3B_j)\Big(~\fiint_{3B_j}|\nabla u(Y)|^2\,dm(Y)\Big)^{\frac12} \lesssim \Big(\dint_{3B_j}\delta(Y)^2|\nabla u(Y)|^2\,\frac{dm(Y)}{m(B_Y)}\Big)^{\frac12}
\end{split} \end{equation}
where we have used the doubling property of $m$ \eqref{mdoubling}, the Poincar\'e inequality \eqref{defPoincare}, and the fact that $\operatorname{rad}(3B_j)\approx\delta(X_j)\approx\delta(Y)$ for each $Y\in 3B_j$. A similar estimate holds for $|u_{B_{j+1}}-u_{3B_j}|$. The combination of \eqref{eq.balls} and \eqref{eq.poincare1} allows us to conclude
\begin{equation}\label{eq.balls2}
1 \lesssim \max_j \dint_{3B_j}\delta(Y)^2|\nabla u(Y)|^2\,\frac{dm(Y)}{m(B_Y)}.
\end{equation}

\subsubsection{A strip decomposition of a wide cone} Recall that $x\in E$ and $S^\ell\in\m O_\ell$ was chosen to contain $x$. The balls $\{B_{j(\ell)}\}$ are the Harnack chain between $\wh B_\ell$ and $\wt B_\ell$ constructed in the beginning of the step. Let us show that there exist $K\geq1$, $\alpha>0$ and an even number $N_K\geq2$ large enough so that for all even $\ell\geq N_K$,
\begin{equation}\label{eq.inclusion}
3B_{j(\ell)}\subset\gamma_{\alpha,\ell}^{\ell(Q)}:=\gamma^{\ell(Q)}_\alpha(x)\cap\big\{Y\in\Omega:\ell(S^\ell)/K\leq\delta(Y)\leq K\ell(S^\ell)\big\}.
\end{equation}
Using the property (c) of the good $\eps_0-$cover, and the fact that $S^\ell\cap S^m\supset\{x\}$ for each $0\leq\ell\leq m$, it is easy to see that  
\begin{equation}\label{eq.smaller}
\ell(S^m)\leq2^{-(m-\ell)}\ell(S^\ell).
\end{equation}
Now, by our constructions we have the chain
\begin{equation}\label{eq.similarities}
\delta(Y) \approx \operatorname{rad}(3B_j)\approx \delta(X_j) \approx\delta(\wt X_\ell)\approx r_\ell = \ell(S^\ell) \qquad\text{ for }Y\in 3B_j,
\end{equation}
and so in particular there exists $K\geq1$ so that $r_\ell/K\leq\delta(Y)\leq Kr_\ell$ for each $Y\in 3B_j$. We fix this $K$. Then, using (\ref{eq.smaller}), we have that $\delta(Y)\leq2^{-\ell}K\ell(Q)$, and so we set $N_K$   even and large enough such that $2^{-N_K}K\leq1$. Hence for all even $\ell\geq N_K$, we have that $\delta(Y)\leq \ell(Q)$. It remains only to find $\alpha$ so that $|Y-x|\leq \alpha \delta(Y)$ for all $Y\in 3B_j$. However, for each $Y\in 3B_j$, armed with (\ref{eq.similarities}) we estimate
\begin{multline}\nonumber
|Y-x|\leq|Y-X_j|+|X_j-\wt X_\ell|+|\wt X_\ell-x_{\wt S^\ell}|+\diam S^\ell\\ \lesssim\operatorname{rad}(3B_j)+\delta(\wt X_{\ell})+\ell(\wt S^\ell)+\ell(S^\ell)\lesssim\delta(Y).
\end{multline}
In summary, $|Y-x|\leq \alpha \delta(Y)$ for some large $\alpha$, as desired. With our choices of $K$, $N_K$, and $\alpha$, (\ref{eq.inclusion}) is proven for all even $\ell\geq N_K$. The combination of \eqref{eq.inclusion} and \eqref{eq.balls2} yields that
\begin{equation}\label{eq.balls3}
1 \lesssim \dint_{\gamma_{\alpha,\ell}^Q}\delta(Y)^2|\nabla u(Y)|^2\,\frac{dm(Y)}{m(B_Y)}.
\end{equation}

\subsubsection{Conclusion of Step III} We are ready to estimate the square function. First, since $k\approx\frac{\log(\zeta^{-1})}{\log(\eps_0^{-1})}\ra\infty$ as $\zeta\ra0$, we consider only $\zeta$ small enough so that $k\geq 4N_K$. Owing to (\ref{eq.smaller}), the strips $\gamma_{\alpha,\ell}^{\ell(Q)}$ have uniformly bounded overlap. Reckon the bounds
\begin{multline}\label{eq.est}
\big|\n A^{Q}_{\alpha}(\delta\nabla u)(x)\big|^2=\dint_{\gamma_{\alpha}^{\ell(Q)}}\delta(Y)^2|\nabla u(Y)|^2\frac{dm(Y)}{m(B_Y)}\\ 
\gtrsim\sum_{\ell=N_K,\, \ell\text{ even}}^k\dint_{\gamma_{\alpha,\ell}^{\ell(Q)}}\delta(Y)^2|\nabla u(Y)|^2\frac{dm(Y)}{m(B_Y)} \gtrsim\sum_{\ell=N_K,\, \ell\text{ even}}^k 1 \gtrsim\frac{k-N_K}2 \approx k,
\end{multline}
where in the second line we used the bounded overlap of the strips, the bound (\ref{eq.balls3}), and the fact that $k\gg N_K$.

\subsection{Step IV: From large square function to $A_{\infty}$} Integrate (\ref{eq.est}) over $x\in E$ with respect to $\sigma$ to see that
\begin{equation}\nonumber
\frac{\log(\zeta^{-1})}{\log(\eps_0^{-1})}\sigma(E)\lesssim k\sigma(E)\lesssim\int_E|\n A_{\alpha}^{Q}(\delta\nabla u)|^2\,d\sigma\leq\int_{Q}|\n A_{\alpha}^{Q}(\delta\nabla u)|^2\,d\sigma\lesssim_{\beta} M\sigma(Q),
\end{equation}
where the last line is a consequence of the assumption \eqref{Main1aloc}. We deduce 
\begin{equation}\label{eq.if}
\frac{\sigma(E)}{\sigma(Q)}\leq C\frac{\log(\eps_0^{-1})}{\log(\zeta^{-1})},\qquad\text{for all Borel }E\subset Q\text{ with }\omega^{X_0}(E)\leq\zeta\omega^{X_0}(Q).
\end{equation}
Given $\xi>0$ and $E\subset Q\in\bb D_{Q_0}$ such that $\omega^{X_0}(E)\leq\zeta\omega^{X_0}(Q)$,  we want to conclude that $\sigma(E)\leq\xi\sigma(Q)$. It is clear that for $\zeta=\zeta(\xi)$ small enough, we achieve the desired result through the estimate (\ref{eq.if}). We have   established that $\omega\in A_{\infty}^{\dyadic}(\sigma,Q_0)$, as desired.\hfill{$\square$}

\section{Proof of Theorem \ref{Main2}} \label{S5}

\begin{lemma} \label{S<N}
Let $(\Omega,m,\mu)$ be PDE friendly. Let $L_0= - \div w\A_0\nabla$ and $L_1 = - \div w\A_1 \nabla$ be two elliptic operators satisfying \eqref{defelliptic} and \eqref{defbounded}, and construct the elliptic measure $\omega_0:= \{\omega_0^X\}_{X\in \Omega}$ and $\omega_1:= \{\omega_1^X\}_{X\in \Omega}$  as in \eqref{defhm}. 

Assume that the weak solutions to $L_1 u = 0$ are the same as the ones of \\$\wh L_1 = -\div w\wh \A_1 \nabla + w\wh \B_1 \cdot \nabla $, and that $\wh \A_1$ still satisfies \eqref{defelliptic}--\eqref{defbounded}. In addition, we require the existence of $K$ such that $\A_0$, $\wh  \A_1$, and $\wh \B_1$ satisfy 
\begin{equation} \label{Main21a}\nonumber
|\wh \A_1- \A_0| \in  {KCM}_{\sup}(\omega_0,K) \quad \text{ and } \quad \delta |\wh \B_1| \in KCM(\omega_0,K).
\end{equation}
Then for any $x\in \partial \Omega$, any $r>0$, any $X \in \Omega \setminus B(x,1000r)$, and any weak solution $u$ to $L_1 u=0$, we have that  
\begin{equation} \label{Main21b}
\int_{\Delta(x,r)} |\n A^r(\delta \nabla u)|^2 d\omega_0^X \leq C(1+K) \int_{\Delta(x,25r)} |N^{2r}(u)|^2 d\omega^X_0,
\end{equation}
where the constants depends only on $n$, the elliptic constants of $\wt \A_0$ and $\wt \A_1$, and the constants in \eqref{mdoubling}, \eqref{Harnack1}, \eqref{dphm1}, and  \eqref{tcp18}.
\end{lemma}

\begin{remark}
The above lemma looks a bit technical, with the introduction of $\wh L_1$. The key observation is that the cases in Theorem \ref{Main2} (multiplicative Carleson perturbation and antisymmetric Carleson perturbation) can be reduced to drift perturbations via the identities \eqref{eq.scalarmult}--\eqref{eq.antisym}, see the proof of Theorem \ref{Main2} below. 

Actually, Lemma \ref{S<N} could be stated without any mention of $L_1$, because the constants in \eqref{Main21b} depends on the properties of $L_0$ and $\wh L_1$, and so only the latter operators matter. 
The only problem lies in the construction of the elliptic measure associated to the $\wh L_1$. In Lemma \ref{S<N}, since $\wh L_1$ has the same solutions as $L_1$, the elliptic measure associated to $\wh L_1$ is the same as $L_1$, hence exists and has the desired properties.

If we had a definition and good properties (the ones presented in Section \ref{S2}) of elliptic measure for (a class of) operators with drifts, for instance by deepening the theory in \cite{dhm}, then we would not really need $L_1$. 
We could only consider two operators with drifts $\wh L_i = -\div w\wh \A_i \nabla + w\wh \B_i \cdot \nabla$, $i\in \{0,1\}$, and their elliptic measures $\omega_i$. And as long as $|\wh \A_1 - \wh \A_0| \in {KCM}_{\sup}(\omega_0)$ and $|\wh \B_1 - \wh \B_0| \in KCM(\omega_0)$, we would have $\omega_1 \in A_\infty(\omega_0)$.  
\end{remark}

\noindent {\em Proof of Theorem \ref{Main2}.}
Since $L_1$ is a (generalized) Carleson perturbation of $L_0$, there exists a function $b$, a matrix $\mathcal C$, and an antisymmetric matrix $\mathcal T$ such that 
\begin{equation} \label{zzeer}
|\mathcal C| \in KCM_{\sup}(\omega_0,K) \quad \text{ and } \frac{\delta |\nabla b|}{b} + \delta w^{-1} |\div(w \mathcal T)| \in KCM(\omega_0,K)
\end{equation}
for some $K>0$, and 
\[\A_1 = b(\A_0 + \mathcal C + \mathcal T).\]
We define 
\[\begin{split}
\wh L_1 
& := -\div( w[\wh \A_0 + \mathcal C] \nabla) - \Big[\div(w\mathcal T) + w \frac{\nabla b}{b}\Big] \cdot \nabla \\
& := - \div( w\wh \A_1 \nabla) - w \wh \B_1 \cdot \nabla.
\end{split}\]
The identities \eqref{eq.scalarmult}--\eqref{eq.antisym} infer that the weak solutions of $L_1$ and $\wh L_1$ are the same. Moreover, \eqref{zzeer} implies that 
\[|\wh \A_1- \A_0| \in  {KCM}_{\sup}(\omega_0,K) \quad \text{ and } \quad \delta |\wh \B_1| \in KCM(\omega_0,K).\]
So we can apply Lemma \ref{S<N} to deduce the bound \eqref{Main21b}. We construct a finitely overlapping covering of $\Delta(x,r)$ by small boundary balls $\{\Delta(x_i,r')\}$ of radius $r'=c_1r/10^6$, where $c_1$ in the constant in the Corkscrew point condition, so that our Corkscrew point $X$ associated to $(x,r)$ stays outside of every $B(x_i,1000r')$. Then, by applying \eqref{Main21b} to every small boundary ball $\Delta(x_i,r')$, we deduce that 
\[\int_{\Delta(x,r)} |\n A^{r'}(\delta \nabla u)|^2 d\omega_0^X \leq C \int_{\Delta(x,2r)} |N^{2r}(u)|^2 d\omega^X_0.
\]
In order to change $\n A^{r'}$ to $\n A^r$ in the above estimate, and hence obtain \eqref{Main2b}, we need to bound the difference $T := |\n A^r(\delta \nabla u)|^2 - |\n A^{r'}(\delta \nabla u)|^2$. We have
\[T(y) = \dint_{W(y,r)} |\delta \nabla u|^2 \frac{dm}{m(B_Y)},\]
where $W(y,r) := \{Y\in \Omega, \, |Y-y| \leq 2\delta(X) \leq 2r, \, r' < \delta(Y)\}$. Notice that all the points $Y\in W(y,r)$ are Corkscrew points associated to $(y,r)$. Therefore, for $Y \in W(y,r)$, we have $\delta(Y) \approx r$ and the doubling property of $m$ infers that 
\[m(B_Y) \approx m(B(y,r) \cap \Omega) \approx m(W(y,r)) \qquad \text{ for } Y \in W(y,r).\]
We conclude that 
\[T(y) \approx  \frac{r^2}{m( W(y,r))} \dint_{W(y,r)} |\nabla u|^2 \, dm.\]
Owning to Caccioppoli's inequality (see for instance Lemma 11.12 in \cite{dfm20}),  one has that $T(y) \lesssim \frac1{m(W^*(y,r))} \iint_{W^*(y,r)} |u|^2 dm$, where $W^*(y,r)$ is a region slightly fatter than $W(y,r)$. From there, it is fairly easy to check that
\[\int_{\Delta(x,r)} T(y) d\omega^{X}_0(y) \lesssim \int_{\Delta(x,r)} \sup_{W^*(y,r)} |u|^2  d\omega^{X}_0(y) \lesssim  \int_{\Delta(x,2r)} |N^{2r}(u)|^2 d\omega^X_0.\]
The theorem follows.
\ep

The rest of the section is devoted to the proof of  Lemma \ref{S<N}. 

\subsection{Step 0: Carleson estimate}
We shall need some preliminary results about the non-tangential maximal function $N$. Note that if one is not interested in the $S<N$ local $L^2$-estimate but only in establishing \eqref{Main1a}, then we could avoid these preliminary estimates and greatly simplify Step 5. But we believe that the $S<N$ estimate is important on its own, and we decided to prove it.

We shall need the untruncated versions of $\n A$ and $N$. We construct the infinite cone $\gamma_\alpha(x) = \{X\in \Omega, \, |X-x| \leq \alpha \delta(X)\}$, and we write $\gamma(x)$ for $\gamma_2(x)$. Then we define, for $f\in L^2_{\loc}(\Omega,m)$ and $x\in \partial \Omega$,
\begin{equation} \label{defAN}\nonumber
\n A(f)(x) := \Big(\dint_{\gamma(x)} |f(x)|^2 \frac{dm(X)}{m(B_X)}\Big)^\frac12 \qquad \text{ and } \qquad N(f)(x) := \sup_{\gamma(x)} |f|.
\end{equation}
We shall also need the variants
\begin{equation} \label{defNvar}\nonumber
\wt N(f)(x) := \sup_{X\in \gamma(x)} \Big(~\fiint_{B_X} |f|^2 \, dm\Big)^\frac12 \qquad \text{ and } \qquad N_{10}(f)(x) := \sup_{\gamma_{10}(x)} |f|.
\end{equation}

Observe that $\wt N(f) \leq N_{10}(f)$, and if we take $2B_X$ instead of $B_X$ in the definition of $\wt N$, the result would still hold. We also have 
\begin{equation} \label{N10<N}
\|N_{10}(f)\|_{L^2(\sigma)} \lesssim \|N(f)\|_{L^2(\sigma)}
\end{equation}
whenever $\sigma$ is doubling on the support of $N_{10}(f)$. The $L^1$ nonlocal result in $\R^n$ can be found in Chapter II, $\S$ 2.5.1 from \cite{ste}, but the proof goes through in our setting without difficulty. The area integral $\n A$, the non-tangential maximal function $N$, and the Carleson measure condition are nicely related via the Carleson inequality. Indeed, if $v\in L^2_{\loc}(\Omega,m)$, $f\in KCM(\sigma,M_f)$ and $\sigma$ is doubling on a neighborhood of the support of $N(v)$, then
\begin{equation} \label{Carleson}
\int_{\partial \Omega} |\n A(fv)|^2 d\sigma \leq C M_f \int_{\partial \Omega} |N(v)|^2 d\sigma,
\end{equation}
where $C$ depends only on the doubling constant of $\sigma$. If $f\in   KCM_{\sup}(\sigma,M_f)$ instead, we can use \eqref{Carleson} with $\wt f(X) = \sup_{B_X} f$ and $\wt v = (\dashint_{B_X} |v^2|dm )^{1/2}$ and obtain the variant
\begin{equation} \label{Carleson2}
\int_{\partial \Omega} |\n A(fv)|^2 d\sigma \leq C M_f \int_{\partial \Omega} |\wt N(v)|^2 d\sigma.
\end{equation}
The proof of \eqref{Carleson} is classical, see for instance \cite[Section II.2.2, Theorem 2]{ste} for the proof on the upper half plane, but which can easily adapted to our setting.

We fix now once for all the rest of this section $x\in \partial \Omega$ and $r>0$.

\subsection{Step 1: Construction of the cut-off function $\Psi$.}  
We choose   a function $\psi \in C^\infty_c(\R)$ that satisfies $0 \leq\psi \leq 1$, $\psi \equiv 1$ on $(-1,1)$, $\psi \equiv 0$ outside $(-2,2)$, and  $|\psi'| \leq 2$. We construct $\Psi = \Psi_{x,r}$ on $\Omega$ as
\begin{equation} \label{defPsi}\nonumber
\Psi(Y) = \psi \Big( \frac{\dist(Y,\Delta(x,r))}{4\delta(Y)} \Big) \psi \Big( \frac{\delta(Y)}{r} \Big)
\end{equation}
and then 
\begin{equation} \label{defPsi2}\nonumber
\Psi_\epsilon(Y) = \Psi(Y) \psi \Big(\frac{\epsilon}{\delta(Y)} \Big).
\end{equation}
Observe that for any $y\in \Delta(x,r)$ and any $Y\in \gamma^r(y)$, we have $\Psi(Y) = 1$. That is, for any $X\in \Omega$, we have
\begin{equation} \label{S<Nb}\nonumber
\int_{\Delta(x,r)} |\n A^r(\delta \nabla u)|^2 d\omega_0^X \leq \int_{\partial \Omega} |\n A(\Psi^2 \delta \nabla u)|^2 d\omega_0^X = \lim_{\epsilon \to 0} \int_{\partial \Omega} |\n A(\Psi_\epsilon^2 \delta \nabla u)|^2 d\omega_0^X.
\end{equation}
Remark also that $\Psi(Y) \neq 0$ means that $\dist(Y,\Delta(x,r)) \leq 8\delta(Y) \leq 16r$, so if $y \in \partial \Omega$ is such that $Y\in \gamma(y)$ and $\Psi(Y) \neq 0$, we necessary have $|y-x| < 21r$. We conclude that
\begin{equation} \label{S<Nc}\nonumber
 \lim_{\epsilon \to 0} \int_{\partial \Omega} |N(\Psi_\epsilon u)|^2 d\omega_0^X =  \int_{\partial \Omega} |N(\Psi u)|^2 d\omega_0^X \leq \int_{\Delta(x,25r)} |N^{2r}(u)|^2 d\omega_0^X.
\end{equation}
As a consequence, \eqref{Main21b} will be proved once we establish that, for $\epsilon >0$, we have
\begin{equation} \label{S<Nd}
\int_{\partial \Omega} |\n A(\Psi_\epsilon^2 \delta \nabla u)|^2 d\omega_0^X \lesssim   \int_{\partial \Omega} |N(\Psi_\epsilon u)|^2 d\omega_0^X.
\end{equation}

\subsection{Step 2: Properties of $\Psi_\epsilon$.} 
In this step we show that $|\nabla \Psi_\epsilon| \in  {KCM}_{\sup}(\omega_0)$. Notice that
\begin{equation} \label{S<Ne}
|\nabla\Psi_\epsilon(Y) |  \lesssim \frac1{\delta(Y) } \1_{E_1 \cup E_2 \cup E_3} \qquad \text{ for } Y \in \Omega,
\end{equation}
where 
\[E_1 := \{Y\in \Omega, \dist(Y,\Delta(x,r))/8 \leq \delta(Y) \leq \dist(Y,\Delta)/4\},\]
\[
E_2 := \{Y\in \Omega, \, r \leq \delta(Y) \leq 2r\},\qquad\text{and}\qquad E_3 := \{Y\in \Omega, \, \epsilon/2 \leq \delta(Y) \leq \epsilon\}.
\]
In addition, if $y\in \partial \Omega$, $Y\in \gamma(y)$, $Y'\in B_{Y}$, and $Y' \in E_1$, then $3\delta(Y)/4 \leq \delta(Y') \leq 5\delta(Y)/4$,
\[\begin{split}
\dist(y,\Delta(x,r)) & \geq \dist(Y',\Delta(x,r)) - |Y'-Y| - |Y-y|  \geq 4\delta(Y') -  \frac{1}4 \delta(Y) - 2\delta(Y)  \geq \frac34 \delta(Y),
\end{split}\]
and
\[\begin{split}
\dist(y,\Delta(x,r)) \leq \dist(Y',\Delta(x,r)) + |Y'-Y| + |Y-y| \leq  13 \delta(Y);
\end{split}\]
that is, for $Y \in \gamma(y)$ such that $B_Y \cap E_1 \neq \varnothing$, 
\begin{equation} \label{S<Nf}
\tfrac1{13} \dist(y,\Delta(x,r)) \leq \delta(Y) \leq\tfrac43\dist(y,\Delta(x,r)).
\end{equation}
We write $\wt{\1_{E_1}}$ for the function $Y \to \sup_{B_Y} \1_{E_1}$, the above estimates proves that $\delta(Y) \approx r_y := \dist(y,\Delta(x,r))$ whenever $Y \in \gamma(y) \cap \supp \, \wt{\1_{E_1}}$. As a consequence, for $y\in \partial \Omega$ and $s>0$, we have that 
\[\begin{split}
|\n A^s(\wt{\1_{E_1}})(y)|^2 \lesssim  \dint_{Y\in \gamma(y), \, \delta(Y) \approx r_y} \frac{dm(Y)}{m(B_Y)} \lesssim 1 
\end{split}\]
because \eqref{mdoubling} implies, for all $Y \in S_y:= \{Y\in \gamma(y), \, \delta(Y) \approx r_y\}$, that $m(B_Y) \approx m(S_y) \approx m(B(y,r_y)\cap \Omega)$. The measure $\omega_0$ does not matter to be able to conclude that $\1_{E_1} \in  {KCM}_{\sup}(\omega_0,M)$, where $M$ depends only on $n$ and the constant in \eqref{mdoubling}.

For $y\in \partial \Omega$, $Y\in \gamma(y)$, $B_Y \cap (E_2 \cup E_3) \neq \varnothing$, we easily deduce from the definition of $E_2$ and $E_3$ that $\delta(Y) \approx r$ or $\delta(Y) \approx \epsilon$. Those estimates are the analogue for $E_2$ and $E_3$ of the bounds \eqref{S<Nf}. With the same arguments as the one used for $E_1$, we obtain that $\1_{E_2\cup E_3}  \in  {KCM}_{\sup}(\omega_0,M)$, hence
\begin{equation} \label{S<Ng}
\1_{E_1 \cup E_2\cup E_3}  \in  {KCM}_{\sup}(\omega_0,M). 
\end{equation}

We combine \eqref{S<Ng} with \eqref{S<Ne} to conclude that
\begin{equation} \label{S<Nh}
|\delta \nabla \Psi_\epsilon|^{1/2} + |\delta  \nabla \Psi_\epsilon|  \in  {KCM}_{\sup}(\omega_0,M) 
\end{equation}
with a constant $M$ that depends only on $n$ and the constant in \eqref{mdoubling}, as desired. Of course, we also have the weaker version 
\begin{equation} \label{S<Nh2}
|\delta \nabla \Psi_\epsilon|^{1/2} + |\delta \nabla \Psi_\epsilon|  \in {KCM}(\omega_0,M). 
\end{equation}

\subsection{Step 3: Introduction of the Green function.} 
The pole $X$ of the elliptic measure $\omega_0$ is chosen in $\Omega \setminus B(x,1000r)$ as in the assumption of the lemma. As an intermediate tool, we shall call $G^*_X$ the weak solution to $(L_0)^*u=0$ in $B(x,500r) \cap \Omega$ that satisfies \eqref{tcp18}. More precisely, we have $\iint_\Omega \A_0\nabla \varphi \cdot \nabla G^*_X \, dm = 0$ for each $\varphi \in C^\infty_c(B(x,500r) \cap \Omega)$, and for $y\in  \Delta(x,25r)$, $s\in (0,2r)$, and any Corkscrew point $Y$ associated to $(y,s)$, the bounds \eqref{tcp18} show that
\begin{equation} \label{S<Nj}
C^{-1} \omega_0^X(\Delta(y,s))  \leq  \frac{m(B(y,s) \cap \Omega)}{s^2} G^*_X(Y) \leq C \omega_0^X(\Delta(y,s)).
\end{equation}

The Green function will be used to replace  the expression with the functional $\n A$ by some integrals over $\Omega$. We claim that, for any $v \in L^2_{\loc}(\Omega)$, we have
\begin{equation} \label{S<Nk}
\int_{\partial \Omega} |\n A(\Psi v)|^2 d\omega_0^X \approx  \dint_\Omega \Psi^2 v^2 \frac{G^*_X}{\delta^2} \, dm.
\end{equation}
Observe that $Y\in \gamma(y)$ implies that $y\in 8\overline{B_Y}\cap\partial\Omega$. As a consequence, Fubini's lemma entails that
\[\int_{\partial \Omega} |\n A(\Psi v)|^2 d\omega_0^X \approx \dint_{\Omega} \Psi^2(Y) v^2(Y)  \frac1{m(B_Y)} \omega_0^X(8B_Y \cap \partial \Omega) dm(Y).\]
Take $Y$ to be such that $\Psi(Y) \neq 0$, and then take $y\in \partial \Omega$ and $s>0$ be such that $s = |y-Y| = \delta(Y)$. The study in Step 1 showed that $y\in \Delta(x,25)$ and $s<2r$, so in particular $X\in \Omega \setminus B(y,2s)$. 
The doubling property of $\omega_0^X$ \eqref{dphm1} shows that $\omega_0^X(8B_Y \cap \partial \Omega) \approx \omega (\Delta(y,s))$, and the doubling property of $m$, given by \eqref{mdoubling}, entails that $m(B_Y) \approx m(B(y,s)\cap \Omega)$. Combined with \eqref{S<Nj}, 
\[ \frac1{m(B_Y)} \omega_0^X(8B_Y \cap \partial \Omega) \approx \frac{G^*_X(Y)}{\delta(Y)^2}.\]
The claim \eqref{S<Nk} follows.

\subsection{Step 4: Bound on the square function.} 
As explained in Step 1, we need to prove \eqref{S<Nd} for any $\epsilon>0$. We define  
\[I = I_\epsilon := \int_{\partial \Omega} |\n A(\Psi_\epsilon^2\delta \nabla u)|^2 d\omega_0^X,\]
which is the quantity that we want to bound. We also set
\begin{equation} \label{defJS<N}\nonumber
J = J_\epsilon := \int_{\partial \Omega} \Big|\wt N\Big(u \Psi_\epsilon^2 \frac{\delta \nabla G^*_X}{G^*_X}\Big) \Big|^2 d\omega^X_0 + \int_{\partial \Omega}  |N(u \Psi_\epsilon  )  |^2 d\omega^X_0.
\end{equation}
If $K$ is the constant in Theorem \ref{Main2},  we claim that, 
\begin{equation} \label{claimS<N}
I \lesssim (1+K)^{1/2} I^{1/2} J^{1/2} + J,
\end{equation}
which self-improves, since $I$ is finite, to $I \lesssim (1+K) J$, or
\begin{equation} \label{S<Nl}
\int_{\partial \Omega} |\n A(\Psi_\epsilon^2\delta \nabla u)|^2 d\omega_0^X \lesssim (1+K) \int_{\partial \Omega} \Big|\wt N\Big(\Psi_\epsilon^2 u \frac{\delta \nabla G^*_X}{G^*_X}\Big) + N(\Psi_\epsilon u)  \Big|^2 d\omega^X_0. 
\end{equation}

Thanks to \eqref{S<Nk}, we have
\begin{equation} \label{S<Nm}
I \approx \dint_\Omega \Psi_\epsilon^4 |\nabla u|^2 G^*_X \, dm.
\end{equation}
Using the ellipticity of $\wh \A_1$, we have
\begin{multline*}
I   \lesssim \dint_\Omega  \wh \A_1 \nabla u \cdot \nabla u  \Psi_\epsilon^4 G^*_X \, dm \\
  = \dint_\Omega \wh \A_1 \nabla u \cdot \nabla [u\Psi_\epsilon^4 G^*_X] \, dm - 4 \dint_\Omega \wh \A_1 \nabla u \cdot \nabla \Psi_\epsilon \, [u \Psi_\epsilon ^3 G^*_X] \, dm  - \dint_\Omega \wh \A_1 \nabla u \cdot \nabla G^*_X \, [u \Psi_\epsilon ^4] \, dm \\
  =: I_1 + I_2 + I_3.
\end{multline*}
We use the fact that $u$ is a weak solution to $L_1$, and thus to $\wh L_1$, to write that  
\[\begin{split}
I_1 & = -\dint_\Omega \wh \B_1\cdot  \nabla u \, [u\Psi_\epsilon^4 G^*_X]  \, dm.
\end{split}\]
We use the Cauchy-Schwarz inequality, \eqref{S<Nm}, and then \eqref{S<Nk} to obtain
\begin{multline*}
I_{1}   \leq \Big(\dint_\Omega  \Psi_\epsilon^4 |\nabla u|^2 G^*_X\, dm \Big)^\frac12   \Big(\dint_\Omega  |\wh \B_1|^2 u^2 \Psi_\epsilon^4 G^*_X\, dm \Big)^\frac12 \\
  \lesssim I^{1/2} \Big(\int_{\partial \Omega}  |\n A(\delta |\wh \B_1| u \Psi^2_\epsilon)|^2  \, d\omega_0^X  \Big)^\frac12 \lesssim I^{1/2} K^{1/2} \Big(\int_{\partial \Omega}  |N(u\Psi_\epsilon)|^2  \, d\omega_0^X  \Big)^\frac12  \lesssim I^{1/2} K^{1/2} J^{1/2},
\end{multline*}
 where the last line is due to the Carleson inequality \eqref{Carleson} and the fact that $\delta|\wh \B_1| \in KCM(\omega_0,K)$.  For $I_2$, the argument is similar, but instead we use the fact that $\wh\A_1$ is bounded, and then the fact that $\nabla \Psi \in KCM(\omega_0,M)$, proved previously in \eqref{S<Nh2}, to get
 \[\begin{split}
I_2 & \lesssim I^{1/2} \Big(\dint_\Omega  |\nabla \Psi_\epsilon|^2 u^2 \Psi_\epsilon^2 G^*_X\, dm \Big)^\frac12 \lesssim I^{1/2} \Big(\int_{\partial \Omega}  |\n A(\delta |\nabla \Psi_\epsilon| u \Psi_\epsilon)|^2  \, d\omega_0^X  \Big)^\frac12 \\
& \lesssim I^{1/2} \Big(\int_{\partial \Omega}  |N(u\Psi_\epsilon)|^2  \, d\omega_0^X  \Big)^\frac12 \lesssim I^{1/2} J^{1/2}.
\end{split}\] 
For the term $I_3$, we replace $\wh \A_1$ by $\A_0$:
 \[ 
I_3   = - \dint_\Omega \A_0 \nabla u \cdot \nabla G^*_X \, [u \Psi_\epsilon ^4] \, dm  - \dint_\Omega (\wh \A_1- \A_0) \nabla u \cdot \nabla G^*_X \, [u \Psi_\epsilon ^4] \, dm  
  = I_{31} + I_{32}. \] 
We deal with $I_{32}$ by invoking the assumption $|\wh \A_1- \A_0| \in  {KCM}_{\sup}(K)$. We have, using the Cauchy-Schwarz inequality, \eqref{S<Nm}, and \eqref{S<Nk} as before, that
 \[
 \begin{split}
I_{32} & \lesssim I^{1/2} \Big(\int_{\partial \Omega}  \Big|\n A\Big( |\wh \A_1- \A_0| u \Psi_\epsilon^2  \frac{\delta \nabla G^*_X}{G^*_X}\Big)\Big|^2  \, d\omega_0^X  \Big)^\frac12 \\
& \lesssim I^{1/2} K^{1/2} \Big(\int_{\partial \Omega} \Big|\wt N\Big(\Psi_\epsilon^2 u \frac{\delta \nabla G^*_X}{G^*_X}\Big)\Big|^2 d\omega^X_0\Big)^{1/2}\lesssim I^{1/2} K^{1/2} J^{1/2}
\end{split}
\] 
by \eqref{Carleson2}, since $|\wh \A_1- \A_0| \in  {KCM}_{\sup}(K)$.
It remains to bound $I_{31}$. We force everything into the first gradient, and we get that
 \[ 
I_{31}   =  - \frac12 \dint_\Omega \wh \A_0 \nabla [u^2 \Psi_\epsilon ^4]  \cdot \nabla G^*_X \, dm + 2 \dint_\Omega \wh \A_0 \nabla \Psi_\epsilon  \cdot \nabla G^*_X \, [u^2 \Psi_\epsilon^3]\, dm \\
  := I_{311} + I_{312}.
\] 
The integral $I_{311}$ is 0. Indeed $G^*_X$ is a weak solution to $(L_0)^*$, and moreover $u^2 \Psi_\epsilon^4$ is a valid test function because it is compactly supported in $\Omega \setminus \{X\}$ and $u^2\Psi^3 \in W^{1,2}(\Omega,m)$ [remember that $u$ is a solution, so $u$ is locally bounded]. As for $I_{312}$, we use the boundedness of $\A_0$ and the inequality $2ab \leq a^2 + b^2$ to infer that
 \[\begin{split}
I_{312} & \lesssim  \dint_\Omega |\nabla \Psi_\epsilon| \Big( \Psi_\epsilon^2 + \Psi_\epsilon^4 \frac{\delta^2 |\nabla G^*_X|^2}{|G^*_X|^2}\Big)  \, u^2 \frac{G^*_X}{\delta} \, dm \\
& \lesssim \int_{\partial \Omega}  |\n A( |\delta \nabla \Psi_\epsilon|^{1/2}  u \Psi_\epsilon)|^2  \, d\omega_0^X + \int_{\partial \Omega}  \Big|\n A\Big( |\delta \nabla \Psi_\epsilon|^{1/2}  u \Psi^2 _\epsilon \frac{\delta \nabla G^*_X}{G^*_X} \Big)\Big|^2  \, d\omega_0^X.
\end{split}\] 
By \eqref{S<Nh} and \eqref{S<Nh2}, that is, by the fact that $|\delta \nabla \Psi_\epsilon|^{1/2} \in KCM(\omega_0,M)$ and $|\delta \nabla \Psi_\epsilon|^{1/2} \in   {KCM}_{\sup}(\omega_0,M)$, and the Carleson inequalities \eqref{Carleson}--\eqref{Carleson2} we conclude that 
\[ I_{312} \lesssim \int_{\partial \Omega} \Big|\wt N\Big(u \Psi_\epsilon^2 \frac{\delta \nabla G^*_X}{G^*_X}\Big) \Big|^2 d\omega^X_0 + \int_{\partial \Omega} |N (u \Psi_\epsilon )  |^2 d\omega^X_0 = J.\]
The claim \eqref{claimS<N} follows.

\subsection{Step 5: A Caccioppoli inequality.} 
From \eqref{S<Nl} and \eqref{S<Nd}, it remains to check that
\begin{equation} \label{S<Nn}
\int_{\partial \Omega} \Big|\wt N\Big(\Psi_\epsilon^2 u \frac{\delta \nabla G^*_X}{G^*_X}\Big) \Big|^2 d\omega^X_0 \lesssim \int_{\partial\Omega}  |N(\Psi_\epsilon u)|^2 d\omega^X_0. 
\end{equation}

We take $Y \in \Omega$ such that $2B_Y \cap \supp \, \Psi \neq 0$, and we observe that $4B_Y$ does not contain $X \in \Omega \setminus B(x,1000r)$. We want to prove the following variant of the Caccioppoli inequality: 
\begin{equation} \label{claimS<N2}
\fiint_{B_Y} \Psi_\epsilon^4 u^2 \frac{\delta^2 |\nabla G^*_X|^2}{|G^*_X|^2} \, dm \lesssim \fiint_{2B_Y} \Psi_\epsilon^2 u^2 \, dm.
\end{equation}
Recall that $G^*_X$ is positive (easy consequence of \eqref{S<Nj}) and a solution to $(L_0)^*u = 0$ on $4B_Y$. The Harnack inequality (Lemma \ref{Harnack}) yields that
\begin{equation} \label{S<No}
G^*_X(Z) \approx G^*_X(Y) \qquad \text{ for } Z \in 2B_Y.
\end{equation}
Using also the property that $\delta \approx \delta(Y)$ on $B_Y$, the claim \eqref{claimS<N2} is equivalent to the estimate
\begin{equation} \label{claimS<N3}\nonumber
\fiint_{B_Y} \Psi_\epsilon^4 u^2 |\nabla G^*_X|^2 \, dm \lesssim \Big(\frac{G^*_X(Y)}{\delta(Y)}\Big)^2 \fiint_{2B_Y} \Psi_\epsilon^2 u^2 \, dm.
\end{equation}
We construct a cut-off function $\Phi = \Phi_Y$, using the smooth function $\psi$ introduced in Step 1, by $\Phi(Z) = \psi\big(\frac{4|Z-Y|}{\delta(Y)}\big)$. 
Note that $\Phi$ is supported in $2B_Y$,  $\Phi \equiv 1$ on $B_Y$, and $|\nabla \Phi| \leq \delta^{-1}(Y)$. So our claim \eqref{claimS<N4} will be proven once we show that 
\begin{equation} \label{claimS<N4}
T:= \dint_{\Omega} \Psi_\epsilon^4 \Phi^2 u^2 |\nabla G^*_X|^2 \, dm \lesssim U := \frac{|G^*_X(Y)|^2}{\delta(Y)^2}  \dint_{2B_Y} \Psi_\epsilon^2 u^2 \, dm.
\end{equation}

We shall prove that $T \lesssim T^{1/2} U^{1/2}$, which self-improves to \eqref{claimS<N4} because $T$ is finite. Using the ellipticity of $\A_0$ given by \eqref{defelliptic}, we have that 
\[\begin{split}
T & \lesssim \dint_\Omega \A_0 \nabla G^*_X \cdot \nabla G^*_X \, [\Psi_\epsilon^4 \Phi^2 u^2] \, dm \\
& = \dint_\Omega \A_0 \nabla [G^*_X \Psi_\epsilon^4 \Phi^2 u^2] \cdot \nabla G^*_X  \, dm 
 - 4\dint_\Omega \A_0 \nabla \Psi_\epsilon \cdot \nabla G^*_X \, [G^*_X \Psi_\epsilon^3 \Phi^2 u^2] \, dm \\
& \quad - 2 \dint_\Omega \A_0 \nabla \Phi \cdot \nabla G^*_X \, [G^*_X \Psi_\epsilon^4 \Phi u^2] \, dm
- 2 \dint_\Omega \A_0 \nabla  u \cdot \nabla G^*_X \, [G^*_X \Psi_\epsilon^4 \Phi^2 u] \, dm \\
& =: T_1 + T_2 + T_3 + T_4.
\end{split}\]
The term $T_1$ equals 0, because $G^*_X$ is a solution to $(L_0)^*$. Using the boundedness of $\A_0$ and the Cauchy-Schwarz inequality, the term $T_3$ is bounded as
\[
T_3 \lesssim T^{1/2} \Big( \dint_\Omega |G^*_X|^2 \Psi_\epsilon^4 |\nabla \Phi|^2 u^2 \, dm \Big)^{1/2} \lesssim  T^{1/2} U^{1/2}
\]
by \eqref{S<No} and   $|\nabla \Phi| \lesssim \delta^{-1}(Y)$. With the same arguments, we treat $T_2$ as follows
\[
T_2 \lesssim T^{1/2}  \Big( \dint_\Omega |G^*_X(Y)|^2 \Psi_\epsilon^2 \Phi^2 |\nabla \Psi_\epsilon|^2 u^2 \, dm \Big)^{1/2} \lesssim  T^{1/2} U^{1/2}
\]
by \eqref{S<No} and \eqref{S<Ne}, i.e. the fact that $|\nabla \Psi_\epsilon| \lesssim \delta^{-1}(Y)$. As for $T_4$, we have
\begin{equation} \label{S<Np}
T_4 \lesssim T^{1/2}  \Big( \dint_\Omega |G^*_X(Y)|^2 \Psi_\epsilon^4 \Phi^2 |\nabla u|^2 \, dm \Big)^{1/2} \lesssim  |G^*_X(Y)| T^{1/2} \Big( \dint_\Omega \Psi_\epsilon^4 \Phi^2 |\nabla u|^2 \, dm \Big)^{1/2}
\end{equation}
If we write $V = \iint_\Omega \Psi_\epsilon^4 \Phi^2 |\nabla u|^2 \, dm$, 
then the same argument as for $T$, using the fact that $u$ is a weak solution to $L_1$ and $|\nabla \Psi_\epsilon| + |\nabla \Phi| \lesssim \delta(Y)$ on $2B_Y$, yields that 
\[
V  \lesssim V^{1/2} \delta^{-1}(Y) \Big(  \dint_{2B_Y} \Psi_\epsilon^2 u^2 \, dm \Big)^{\frac12} = |G^*_X(Y) |^{-1} V^{1/2} U^{1/2} ,
\]
which self-improves to $ |G^*_X(Y)|^2 V \lesssim U$. Using the estimate in \eqref{S<Np}, we obtain that $T_4 \lesssim T^{1/2} U^{1/2}$. The claim \eqref{claimS<N4} follows, and hence so does \eqref{claimS<N2}.
 
The inequality \eqref{claimS<N2} entails the pointwise bound
\[ \wt N\Big(\Psi_\epsilon^2 u \frac{\delta \nabla G^*_X}{G^*_X}\Big)  \lesssim N_{10} (\Psi_\epsilon u ),\]
where $N_{10}(v)(y) := \sup_{\gamma_{10}(y)} |v|$, and $\gamma_{10}(y)$ is the cone with vertex at $y \in \partial \Omega$ with a bigger aperture than $\gamma(y)$ so that  $\gamma_{10}(y) \supset \bigcup_{Y \in \gamma(y)} 2B_Y$. The estimate \eqref{S<Nn} comes then from classical fact that $\|N_{10}(v)\|_{L^2} \lesssim \|N(v)\|_{L^2}$, see \eqref{N10<N}. If we want to avoid this latter estimate, we can also define $N$ using cones with bigger apertures than the ones of $\wt N$, and all our proofs are then identical.\ep

\hypersetup{linkcolor=toc}

\bibliography{refs} 
\bibliographystyle{alpha-sort-max} 

\end{document}